\documentclass[a4paper, 12pt]{article}
\usepackage{enumerate, theorem}
\usepackage{amsmath, amsfonts, amssymb}
\usepackage[height=22.5cm, width=15cm]{geometry}
\usepackage[all]{xy}
\usepackage{mathrsfs, yfonts}

\DeclareMathOperator{\C}{\mathbb{C}}

\newcommand{\A}{\tilde{\mathcal{A}}}

\newcommand{\parag}[1]{\paragraph{\sc{#1.}}}

\newtheorem{thm}{Theorem}[subsection]
\newtheorem{defn}[thm]{Definition}
\newtheorem{cor}[thm]{Corollary}
\newtheorem{prop}[thm]{Proposition}
\newtheorem{lemma}[thm]{Lemma}

\begin{document}

\title{Asymptotics of a vanishing period : \\
General existence theorem and basic properties  of frescos.}

\author{Daniel Barlet\footnote{Barlet Daniel, Institut Elie Cartan UMR 7502  \newline
Universit\'e de Lorraine, CNRS, INRIA  et  Institut Universitaire de France, \newline
BP 239 - F - 54506 Vandoeuvre-l\`es-Nancy Cedex.France. \newline
e-mail : Daniel.Barlet@iecn.u-nancy.fr}.}

\date{5/1/12.}

\maketitle

\section*{Abstract}
In this paper we introduce the word "fresco" to denote a \  $[\lambda]-$primitive monogenic geometric (a,b)-module. The study of this "basic object" (generalized Brieskorn module with one generator) which corresponds to the minimal filtered (regular) differential equation satisfied by a relative de Rham cohomology class, began in [B.09] where the first structure theorems are proved. Then in [B.10] we introduced the notion of theme which corresponds in the \ $[\lambda]-$primitive case  to frescos having a unique Jordan-H{\"o}lder sequence. Themes correspond  to asymptotic expansion of a given vanishing period, so to the image of a fresco in the module of  asymptotic expansions. For a fixed relative de Rham cohomology class (for instance given by  a smooth differential form $d-$closed and $df-$closed) each choice of a vanishing cycle in the spectral eigenspace of the monodromy for the eigenvalue \ $exp(2i\pi.\lambda)$ \ produces a \ $[\lambda]-$primitive theme, which is a quotient of the fresco associated to the given relative de Rham class  itself. \\
The first part of this paper shows that, for any \ $[\lambda]-$primitive fresco there exists an unique Jordan-H{\"o}lder sequence (called the principal J-H. sequence) with corresponding quotients giving the opposite of  the roots of the Bernstein polynomial in a non decreasing order. Then we introduce and study  the semi-simple part of a given fresco and we  characterize the semi-simplicity of a fresco  by the fact  for any given order of the roots of its  Bernstein polynomial we may find a J-H. sequence making them appear with this order.
Then, using the parameter associated to a rank \ $2$ \ \ $[\lambda]-$primitive theme, we introduce inductiveley a numerical invariant, that we call the \ $\alpha-$invariant, which depends polynomially on the isomorphism class of a  fresco (in a sens which has to be defined) and which allows to give an inductive way to produce a sub-quotient rank \ $2$ \ theme of a given \ $[\lambda]-$primitive fresco assuming non semi-simplicity.\\
In the last section  we prove a general existence result which naturally  associate  a fresco to any relative de Rham cohomology class of a proper holomorphic  function of a complex manifold onto a disc. This is, of course, the motivation for the study of frescos.

\parag{AMS Classification} 32 S 25, 32 S 40, 32 S 50.

\parag{Key words} Fresco,  theme,  (a,b)-module, asymptotic expansion, vanishing period, Gauss-Manin connection, filtered differential equation.

\tableofcontents

\section*{Introduction}

Let \ $f : X \to D$ \ be an holomorphic function on a connected  complex manifold. Assume that \ $\{df = 0 \} \subset \{ f = 0 \}: = X_0 $. We consider \ $X$ \ as a degenerating family of complex manifolds parametrized by \ $D^* : = D \setminus \{0\}$ \ with a singular member  \ $X_0$ \ at the origin of \ $D$. Let \ $\omega$ \ be a smooth \ $(p+1)-$differential form on \ $X$ \ satisfying \ $d\omega = 0 = df\wedge \omega$. Then in many interesting cases (see for instance [B.II] , [B.III] for the case of a function with 1-dimensional singular set and the section 4  for the general proper case) the relative family of de Rham cohomology classes induced on the fibers \ $(X_s)_{ s \in D^*}$ \ of \ $f$ \ by  \ $\omega\big/df$ \ is solution of a minimal filtered differential equation defined from the Gauss-Manin connection of \ $f$. This object, called a {\bf fresco} is a monogenic regular (a,b)-module satisfying an extra condition, called "geometric", which encodes simultaneously the regularity at \ $0$ \ of the Gauss-Manin connection, the monodromy theorem and B. Malgrange's positivity theorem.\\
We study the structure of such an object  in order to determine the possible quotient themes of a given fresco. Such a theme corresponds to a possible asymptotic expansion of vanishing periods constructed from \ $\omega$ \ by choosing a vanishing cycle \ $\gamma \in H_p(X_{s_0},\C)$ \ and definig
$$ F_{\gamma}(s) : = \int_{\gamma_s} \ \omega\big/df  $$
where \ $\gamma_s$ \ is the (multivalued) horizontal family of cycles defined from \ $\gamma$ \  in the fibers of \ $f$ \ (see [M.74]).\\

Let me describe the content of this article.\\
After some easy preliminaries, we prove in section 1 that a \ $[\lambda]-$primitive fresco \ $E$ \  admits an {\em unique} Jordan-H{\"o}lder sequence, called the principal J-H. sequence, in which the opposite of the roots of the Bernstein polynomial of \ $E$ \ appears in a non decreasing order. This uniqueness result is important because it implies, for instance, that the isomorphism classes of each quotient of terms of the principal J-H. sequence only depends on the isomorphism class of \ $E$.

\smallskip

In section 2 we define and study semi-simple regular (a,b)-modules and the corresponding semi-simple filtration. In the case of a  \ $[\lambda]-$primitive fresco we prove that semi-simplicity is characterized by the fact that we may find a J-H. sequence in which the opposite of the roots of the Bernstein polynomial appears in strictely decreasing order (but also in any given order).

\smallskip

In section 3 we answer to the following question : How to recognize from the principal J-H. sequence of a \ $[\lambda]-$primitive fresco if it is semi-simple. Of course the answer is obvious if some rank \ $2$ \ sub-quotient theme appears from this sequence. But when it is not the case, such a rank \ $2$ \ sub-quotient theme may appear after commuting some terms in the J-H. sequence. The simplest example is when
$$ E : = \A\big/\A.(a - \lambda_1.b)(1 + \alpha.b^{p_1+p_2})^{-1}.(a - \lambda_2.b).(a - \lambda_3.b) $$
with \ $\lambda_{i+1} = \lambda_i + p_i -1$ \ for \ $i = 1,2$ \ with \ $p_i \in \mathbb{N}^*,  i = 1,2$ \ and \ $\alpha \in \C^*$. Using the identity in \ $\A$
 $$(a - \lambda_2.b).(a - \lambda_3.b) = (a - (\lambda_3+1).b).(a - (\lambda_2-1).b) $$
 it is easy to see that \ $(a - (\lambda_2-1).b).[1]$ \ generates a (normal) rank \ $2$ \ theme in \ $E$, because we assume \ $\alpha \not= 0$.
 
 \smallskip
 
We solve this question introducing for a \ $[\lambda]-$primitive rank \ $k \geq 2$ \  fresco \ $E$ \ such \ $F_{k-1}$ \ and \ $E\big/F_1$ \ are semi-simple, where \ $(F_j)_{j \in [1,k]}$ \ is the principal J-H. sequence of \ $E$, the \ $\alpha-$invariant \ $\alpha(E)$. This complex number only depends on the isomorphism class of \ $E$ \ and is zero if and only \ $E$ \ is semi-simple. We also prove that when \ $\alpha(E)$ \ is not zero it gives the parameter of any normal rank \ $2$ \ sub-theme of \ $E$.
 
 \smallskip
 
  In section 4 we give, using tools introduced in [B.II], a general existence theorem for the fresco associated to a relative de Rham cohomology class in the geometric situation described in the begining of this introduction, assuming the function \ $f$ \ proper.

\section{Preliminaries.}

\subsection{Definitions and characterization as \ $\A-$modules.}

We are interested in "standard " formal asymptotic expansions of the following type
$$ \sum_{q = 1}^p \sum_{j=0}^N \ \C[[s]].s^{\lambda_q-1}.(Log\, s)^j $$
where \ $\lambda_1, \dots, \lambda_q$ \ are positive rational numbers, and in fact in vector valued such expansions. The two basic operations on such expansions are
\begin{itemize}
\item the multiplication by \ $s$ \ that  we shall denote \ $a$,
\item and the primitive in \ $s$ \ without constant that we shall denote \ $b$.
\end{itemize}
This leads to consider on the set of such expansions a left module structure on the \ $\C-$algebra 
$$ \A : = \{\  \sum_{\nu = 0}^{\infty} \ P_{\nu}(a).b^{\nu},\quad {\rm where} \quad P_{\nu} \in \C[x] \  \}  $$
defined by  the following conditions
\begin{itemize}
\item The commutation relation \ $a.b - b.a = b^2$ \ which is the translation of the Leibnitz rule ;
\item The continuity for the \ $b-$adic filtration of \ $\A$ \ of the left and right multiplications by \ $a$.
\end{itemize}
Define now for \ $\lambda$ \ in \ $\mathbb{Q}\, \cap\, ]0,1]$ \ and \ $N \in \mathbb{N}$ \  the left \ $\A-$module
$$ \Xi^{(N)}_{\lambda} : = \oplus_{j=0}^N \ \C[[s]].s^{\lambda-1}.(Log\, s)^j  = \oplus_{j=0}^N \ \C[[a]].s^{\lambda-1}.(Log\, s)^j = \oplus_{j=0}^N \ \C[[b]].s^{\lambda-1}.(Log\, s)^j .$$
Of course we let \ $a$ \ and \ $b$ \ act on \ $\Xi_{\lambda}^{(N)}$ \ as explained above.\\
Define also, when \ $\Lambda$ \ is a finite  subset in \ $\mathbb{Q}\, \cap\, ]0,1]$, the \ $\A-$module
$$ \Xi_{\Lambda}^{(N)} : = \oplus_{\lambda \in \Lambda}\  \Xi_{\lambda}^{(N)} .$$
More generally, if \ $V$ \ is a finite dimensional complex vector space we shall put a structure of left \ $\A-$module on \ $\Xi_{\Lambda}^{(N)}\otimes_{\C} V$ \ with the following rules :
$$ a.(\varphi \otimes v) = (a.\varphi)\otimes v \quad {\rm and} \quad b.(\varphi \otimes v) = (b.\varphi)\otimes v$$
for any \ $\varphi \in \Xi_{\Lambda}^{(N)}$ \ and any \ $v \in V$. It will be convenient to denote  \ $\Xi_{\lambda}$ \ the \ $\A-$module \ $\sum_{N \in \mathbb{N}}\ \Xi^{(N)}_{\lambda}$. \\

\begin{defn}\label{fresco}
A \ $\A-$module is call a {\bf fresco} when it is isomorphic to a submodule \ $\A.\varphi \subset \Xi_{\Lambda}^{(N)}\otimes V$ \ where \ $\varphi$ \ is any element in \ $\Xi_{\Lambda}^{(N)}\otimes V$, for some choice of \ $\Lambda, N$ \ and \ $V$ \ as above.\\
A fresco is a {\bf theme} when we may choose \ $V : = \C$ in the preceeding choice.
\end{defn}

\bigskip

Now the characterization of frescos among all left \ $\A-$modules is not so obvious. The following theorem is proved in [B.09].

\begin{thm}\label{charact. frescos}
A left \ $\A-$module \ $E$ \  is a fresco if and only if it is a geometric (a,b)-module which is generated (as a \ $\A-$module) by one element. Moreover the annihilator in \ $\A$ \ of  any generator of \ $E$ \ is a left ideal of the form \ $\A.P$ \ where \ $P$ \ may be written as follows
$$ P = (a - \lambda_1.b).S_1^{-1}.(a - \lambda_2.b).S_2^{-1} \dots (a - \lambda_k.b).S_k^{-1}, \quad k : =  \dim_{\C}(E\big/b.E) $$
where \ $\lambda_j$ \ are rational numbers such that \ $\lambda_j + j > k$ \ for \ $j \in [1,k]$ \ and where \ $S_1, \dots , S_k$ \ are invertible elements in the sub-algebra \ $\C[[b]]$ \ of \ $\A$.\\
Conversely, for such a \ $P \in \A$ \ the left \ $\A-$module \ $E : = \A\big/\A.P$ \ is a fresco and it is a free rank \ $k$ \ module on \ $\C[[b]]$. 
\end{thm}

Let me recall briefly for the convenience of the reader the definitions of the notions involved in the previous statement.
\begin{itemize}
\item A {\bf (a,b)-module} \ $E$ \ is a free finite rank \ $\C[[b]]$ \ module endowed with an \ $\C-$linear endomorphism \ $a$ \ such that \ $a.S = S.a + b^2.S'$ \ for \ $S \in \C[[b]]$;  or, in an equivalent way, a \ $\A-$module which is free and finite type over the subalgebra \ $\C[[b]] \subset \A$. \\
 It has a {\bf simple pole} when \ $a$ \ satisfies \ $a.E \subset b.E$. In this case the {\bf Bernstein polynomial} \ $B_E$ \  of \ $E$ \ is defined as the minimal polynomial of \ $-b^{-1}.a$ \ acting on \ $E\big/b.E$.
\item A (a,b)-module \ $E$ \  is {\bf regular}  when it may be embedded in a simple pole (a,b)-module. In this case there is a minimal such embedding which is the inclusion of \ $E$ \ in its saturation \ $E^{\sharp}$ \ by \ $b^{-1}.a$. The {\bf Bernstein polynomial} of a regular (a,b)-module is, by definition, the Bernstein polynomial of its saturation \ $E^{\sharp}$.
\item A regular (a,b)-module \ $E$ \ is called {\bf geometric} when all roots of its Bernstein polynomial are rational and strictly negative.
\end{itemize}

Note that the formal completion in \ $b$ \ of the Brieskorn module of a function with an isolated singularity is a geometric (a,b)-module. The last section of this article shows that this structure appears in a rather systematic way in the study of the Gauss-Manin connection of a proper  holomorphic function on a complex manifold.

\begin{defn}\label{normal}A submodule \ $F $ \ in \ $E$ \ is {\bf normal} when \ $F \cap b.E = b.F$.
\end{defn}

For any sub-module \ $F$ \ of a (a,b)-module \ $E$ \ there exists a minimal normal sub-module \ $\tilde{F}$ \ of \ $E$ \ containing \ $F$. We shall call it the {\bf normalization} of \ $F$. It is easy to see that \ $\tilde{F}$ \ is the pull-back by the quotient map \ $E \to E\big/F$ \ of the \ $b-$torsion of \ $E\big/F$.

\smallskip

 When \ $F$ \ is normal the quotient \ $E\big/F$ \ is again a (a,b)-module. Note that a sub-module of a regular (resp. geometric) (a,b)-module is regular (resp. geometric) and when \ $F$ \ is normal \ $E\big/F$ \ is also regular (resp. geometric).

\begin{lemma}
Let \ $E$ \ be a fresco and \ $F$ \ be a normal sub-module in \ $E$. Then \ $F$ \ is a fresco and also the quotient \ $E\big/F$.
\end{lemma}

\parag{proof} The only point to prove, as we already know that \ $F$ \ and \ $E\big/F$ \ are geometric (a,b)-modules thanks to [B.09], is the fact that \ $F$ \ and \ $E\big/F$ \ are generated as \ $\A-$module by one element. This is obvious for \ $E\big/F$, but not for \ $F$. We shall use the theorem \ref{charact. frescos} for \ $E\big/F$ \ to prove that \ $F$ \ is generated by one element. Let \ $e$ \ be a generator of \ $E$ \ and let \ $P$ \ as in the theorem \ref{charact. frescos} which generates the annihilator ideal of the image of \ $e$ \ in \ $E\big/F$. Then \ $P.e$ \ is in \ $F$. We shall prove that \ $P.e$ generates \ $F$ \ as a \ $\A-$module. Let \ $y$ \ be an element in \ $F$ \ and write \ $y = u.e$ \ where \ $u$ \ is in \ $\A$. As \ $P$ \ is, up to an invertible element in \ $\C[[b]]$, a monic polynomial in \ $a$ \ with coefficients in \ $\C[[b]]$, we may write \ $u = Q.P + R$ \ where \ $Q$ \ and \ $R$ \ are in \ $\A$ \ and \ $R$ \ is a polynomial in \ $a$ \ with coefficient in \ $\C[[b]]$ \ of degree \ $r < deg(P) = rank(E\big/F)$. Now, as \ $y$ \ is in \ $F$,  the image in \ $E\big/F$ \ of \ $u.e$ \ is \ $0$, and this implies that \ $R$ \ annihilates the image of \ $e$ \ in \ $E\big/F$. So \ $R$ \ lies in \ $\A.P$ \ and so \ $R = 0$. Then we have \ $u = Q.P$ \ and \ $y = Q.P.e $ \ proving our claim. $\hfill \blacksquare$\\

In the case of a fresco the Bernstein polynomial is more easy to describe, thanks to the following proposition proved in [B.09].

\begin{prop}\label{Bernst. fresco}
Soit \ $E = \A\big/\A.P$ \ be a rank \ $k$ \ fresco as described in the previous theorem. The Bernstein polynomial of \ $E$ \ is the {\em characteristic polynomial} of \ $-b^{-1}.a$ \ acting on \ $E^{\sharp}\big/b.E^{\sharp}$. And  the {\bf Bernstein element}  \ $P_E$ \ of \ $E$, which is the element in \ $\A$ \ defined by the Bernstein polynomial \ $B_E$ \ of \ $E$ \ by the following formula
$$ P_E : = (-b)^k.B_E(-b^{-1}.a) $$
is equal to \ $(a - \lambda_1.b) \dots (a - \lambda_k.b)$ \ for any such choice of presentation of \ $E$.
\end{prop}

As an easy consequence, the \ $k$ \  roots of the Bernstein polynomial of \ $E$ \ are the opposite of the numbers \ $\lambda_1+1-k, \dots, \lambda_k+k-k$. So the Bernstein polynomial is readable on the element \ $P$ : the initial form of \ $P$ \ in (a,b) is the Bernstein element\ $P_E$ \ of \ $E$. 

\parag{Remark} If we have an exact sequence of (a,b)-modules
$$ 0 \to F \to E \to G \to 0 $$
where \ $E$ \ is a fresco, then \ $F$ \ and \ $G$ \ are frescos and the Bernstein elements satisfy the equality \ $P_E = P_F.P_G$ \ in the algebra \ $\A$ \ (see [B.09]  proposition 3.4.4.). $\hfill \square$

\parag{example} The \ $\A-$module \ $\Xi_{\lambda}^{(N)}$ \ is a simple pole (a,b) with rank \ $N+1$. Its Bernstein polynomial is equal to \ $(x+\lambda)^{N+1}$.\\
The theme \ $\A.\varphi \subset \Xi^{(N)}_{\lambda}$ \ where \ $\varphi = s^{\lambda-1}.(Log\,s)^N$ \ has rank \ $N+1$ \ and Bernstein element \ $(a - (\lambda+N).b)(a - (\lambda+N-1).b) \dots (a - \lambda.b)$. Its saturation is \ $\Xi^{(N)}_{\lambda}$.\\

The dual of a regular (a,b)-module is regular, but duality does not preserve the property of being geometric because duality change the sign of the roots of the Bernstein polynomial. As duality preserves regularity (see [B.95]), to find again a geometric (a,b)-module it is sufficient to make the tensor product  by a rank \ $1$ \ (a,b)-module \ $E_{\delta}$ \ for \ $\delta$ \ a large enough rational number\footnote{to tensor by \ $E_{\delta}$ \ is the same that to replace \ $a$ \ by \ $a + \delta.b$ ; see [B.I] for the general  definition of the tensor product.}. The next lemma states that this "twist duality" preserves the notion of fresco.

\begin{lemma}\label{twisted duality}
Let \ $E$ \ be a fresco and let \ $\delta$ \ be a rational number such that \ $E^*\otimes E_{\delta}$ \ is geometric. Then \ $E^*\otimes E_{\delta}$ \ is a fresco. 
\end{lemma}

The proof is obvious. $\hfill \blacksquare$\\

The following definition is useful when we want to consider only the part of asymptotic expansions corresponding to prescribe eigenvalues of the monodromy.

\begin{defn}\label{primitive}
Let \ $\Lambda$ \ be a subset of \ $\mathbb{Q} \, \cap \, ]0,1]$. We say that a regular (a,b)-module \ $E$ \ is {\bf \ $[\Lambda]-$primitive} when all roots of its Bernstein polynomial are in \ $- \Lambda + \mathbb{Z}$.
\end{defn}

The following easy proposition is proved in [B.09]

\begin{prop}\label{primitive part}
Let \ $E$ \ be a regular (a,b)-module and \ $\Lambda$ \ a subset of \ $\mathbb{Q} \, \cap \, ]0,1]$. Then there exists a maximal submodule \ $E_{[\Lambda]}$ \ in \ $E$ \ which is \ $[\Lambda]-$primitive. Moreover the quotient \ $E\big/E_{[\Lambda]}$ \ is a \ $[\Lambda^c]-$primitive (a,b)-module, where we denote \  $\Lambda^c : = \mathbb{Q} \, \cap \, ]0,1] \setminus \Lambda$.
\end{prop}

We shall mainly consider the case where \ $\Lambda$ \ is a single element. Note that  the \ $[\lambda]-$primitive part of an (a,b)-module  \ $E \subset \Xi^{(N)}_{\Lambda}\otimes V$ \ corresponds to its intersection  with \ $\Xi^{(N)}_{\lambda}\otimes V$.

\subsection{The principal Jordan-H{\"o}lder sequence.}

The classification of rank \ $1$ \ regular (a,b)-module is very simple : each isomorphy class is given by a complex number and to \ $\lambda \in \C$ \ corresponds the isomorphy class of \ $E_{\lambda} : = \A\big/\A.(a - \lambda.b)$. Then {\bf Jordan-H{\"o}lder sequence} for a regular (a,b)-module \ $E$ \ is a sequence 
$$ 0 = F_0 \subset F_1 \subset \dots \subset F_k = E  $$
of normal sub-modules such that for each \ $j \in [1,k]$ \ the quotient \ $F_j\big/F_{j-1}$ \ has rank \ $1$. Then to each J-H. sequence we may associate an ordered sequence of  complex numbers \ $\lambda_1, \dots, \lambda_k$ \ such that \ $F_j\big/F_{j-1} \simeq E_{\lambda_j}$.

\parag{example}
A regular rank 1 (a,b)-module is a fresco if and only if it is isomorphic to \ $E_{\lambda}$ \ for some \ $\lambda \in \mathbb{Q}^{+*}$. All rank 1 frescos are themes.
The classification of rank 2 regular (a,b)-modules given in [B.93] gives the list of \ $[\lambda]-$primitive rank 2 frescos  which is the following, where \ $\lambda_1 > 1$ \ is a rational number : 
\begin{align*}
& E =  E \simeq \A\big/\A.(a -\lambda_1.b).(a - (\lambda_1-1).b)  \tag{1}\\
& E \simeq \A\big/\A.(a -\lambda_1.b).(1+ \alpha.b^p)^{-1}.(a - (\lambda_1 + p-1).b) \tag{2}
\end{align*}
where \ $p \in \mathbb{N} \setminus \{0\}$ \ and \ $\alpha \in \C$.\\
The themes in this list are these in \ $(1)$ \ and these in \ $(2)$ \ with \ $\alpha \not= 0 $. For a \ $[\lambda]-$primitive fresco  in case \ $(2)$ \ the number \ $\alpha$ \ will be called the {\bf parameter} of the rank \ $2$ \ fresco. By convention we shall define \ $\alpha : = 1$ \ in the case \ $(1)$.$\hfill \square$ \\

The following existence result is proved in [B.93]

\begin{prop}\label{J-H exist}
For any regular (a,b)-module \ $E$ \ of rank \ $k$ \ there exists a J-H. sequence. The numbers \ $exp(2i\pi.\lambda_j)$ \ are independant of the J-H. sequence, up to permutation. Moreover the number \ $\mu(E) : = \sum_{j=1}^k \ \lambda_j $ \ is also independent of the choice of the J-H. sequence of \ $E$.
\end{prop}

\parag{Exercice} Let \ $E$ \ be a regular (a,b)-module and  \ $E' \subset E$ \ be a sub-(a,b)-module with the same rank than \ $E$. Show that \ $E'$ \ has finite \ $\C-$codimension in \ $E$ \ given by
$$ \dim_{\C}E\big/E' = \mu(E') - \mu(E) .$$
{\em hint} : make an induction on the rank of \ $E$.$\hfill \square$\\ 

For a \ $[\lambda]-$primitive fresco a more precise result is proved in [B.09]

\begin{prop}
For any J-H. sequence of a rank \ $k$  \ $[\lambda]-$primitive fresco \ $E$ \  the numbers \ $\lambda_j+j-k, j \in [1,k]$ \ are the opposite of the  roots (with multiplicities) of the Bernstein polynomial of \ $E$. So, up to a permutation, they are independant of the choice of the J-H. sequence. Moreover, there always exists a J-H. sequence such that the associated  sequence \ $\lambda_j+j$ \ is non decreasing.
\end{prop}

For a \ $[\lambda]-$primitive theme the situation is extremely rigid (see [B.09]) : 

\begin{prop}\label{theme J-H.}
Let \ $E$ \ a rank \ $k$ \ $[\lambda]-$primitive theme. Then, for each \ $j$ \ in \  $ [0,k]$,  there exists an unique normal rank \ $j$ \ submodule \ $F_j$. The numbers associated to the unique J-H. sequence satisfy the condition that \ $\lambda_j+j, j \in [1,k]$ \ is an non decreasing  sequence.
\end{prop}

\begin{defn}\label{Principal J-H.}
Let \ $E$ \ be a \ $[\lambda]-$primitive fresco of rank \ $k$ \ and let  
 $$0 = F_0 \subset F_1 \subset \dots \subset F_k = E$$
 be a J-H. sequence of \ $E$. Then for each \ $j \in [1,k]$ \ we have \ $F_j\big/F_{j-1} \simeq E_{\lambda_j}$, where \ $\lambda_1, \dots, \lambda_k$ \ are in \ $\lambda + \mathbb{N}$.
We shall say that such a J-H. sequence is {\bf principal} when the sequence \  $[1,k]\ni j \mapsto \lambda_j + j$ \ is non decreasing.
\end{defn}

\begin{prop}\label{Uniqueness}
Let \ $E$ \ be a \ $[\lambda]-$primitive fresco. Then its principal J-H. sequence is unique.
\end{prop}

 We shall prove the uniqueness by induction on the rank \ $k$ \ of \ $E$. 

We begin by  the case of rank 2.

\begin{lemma}\label{petit}
Let \ $E$ \ be a rank 2 \ $[\lambda]-$primitive fresco and let \ $\lambda_1,\lambda_2$ \ the numbers corresponding to a principal J-H. sequence of \ $E$ \ (so  \ $\lambda_1+1 \leq \lambda_2+2$). Then the normal rank 1 submodule of \ $E$ \ isomorphic to \ $E_{\lambda_1}$ \ is unique.
\end{lemma}

\parag{Proof} The case \ $\lambda_1+ 1 = \lambda_2+2$ \ is obvious because then \ $E$ \ is a \ $[\lambda]-$primitive theme (see [B.10] corollary 2.1.7). So we may assume that \ $\lambda_2 = \lambda_1 + p_1-1$ \ with \ $p_1 \geq 1$ \ and that  \ $E$ \ is the quotient \ $E \simeq \A\big/\A.(a - \lambda_1.b).(a - \lambda_2.b)$ \ (see the classification of rank 2 frescos in 2.2), because the result is clear when \ $E$ \ is a theme.
We shall use the \ $\C[[b]]-$basis \ $e_1, e_2$ \  of \ $E$ \ where \ $a$ \ is defined by the relations
$$ (a - \lambda_2.b).e_2 = e_1 \quad (a - \lambda_1.b).e_1 = 0. $$
This basis comes from the isomorphism \ $E \simeq \A\big/\A.(a - \lambda_1.b).(a - \lambda_2.b)$ \  deduced from the classification of rank 2 frescos with \ $e_2 = [1]$ \ and \ $e_1 = (a - \lambda_2.b).e_2$.\\
Let look for \ $x : = U.e_2 + V.e_1$ \ such that \ $(a - \lambda_1.b).x = 0$. Then we obtain
$$ b^2.U'.e_2 + U.(a - \lambda_2.b).e_2 + (\lambda_2-\lambda_1).b.U.e_2 + b^2.V'.e_1 = 0$$
which is equivalent to the two equations :
$$ b^2.U' + (p_1-1).b.U = 0 \quad {\rm and} \quad U + b^2.V' = 0 .$$
The first equation gives \ $U = 0$ \ for \ $p_1 \geq 2$ \ and \ $U \in \C$ \ for \ $p_1 = 1$. As the second equation implies \ $U(0) = 0$, in all cases \ $U = 0$ \ and \ $V \in \C$. So the solutions are in \ $\C.e_1$. $\hfill \blacksquare$\\

Remark that in the previous lemma, if we assume \ $p_1\geq 1$ \ and \ $E$ \ is not a theme, it may exist infinitely many different normal (rank 1) submodules isomorphic to \ $E_{\lambda_2+1}$. But then, \ $\lambda_2+2 > \lambda_1+1$. This happens in the following example.

\parag{Example} Let \ $E : = \A\big/\A.P$ \ with
$$ P : = (a - \lambda_1.b).(a - \lambda_2.b)(1 + \alpha.b^{p_2})^{-1}.(a - \lambda_3.b) $$
where \ $\lambda_1 > 2$ \ is rational, $p_1$ \ and \ $p_2$ \ are in \ $\mathbb{N}^*$, $\lambda_2 = \lambda_1 + p_1 -1, \lambda_3 = \lambda_2 + p_2 - 1$ \ and \ $\alpha$ \ is in \ $\C^*$.\\
Then using the identity (see the commuting lemma in [B.09])
\begin{align*}
&   (a - \lambda_1.b).(a - \lambda_2.b) = 
 U^{-1}.(a - (\lambda_2+1).b).U^2.(a - (\lambda_1-1).b).U^{-1} 
\end{align*}
with \ $U : = 1 + \rho.b^{p_1}, \rho \in \C^*$, it is easy to see that the rank \ $3$ \ fresco \ $E$ \ admits the rank \ $2$ \ theme \ $T$ \ with fundamental invariants \ $(\lambda_1-1, \lambda_3)$ \ and parameter \ $\rho.\alpha$.\\
Remark that if we use the previous identity with \ $\rho = 0$ \ (so \ $U = 1$), then using the  identity
\begin{align*}
&  (a -(\lambda_1-1).b).(1 + \alpha.b^{p_2})^{-1}.(a - \lambda_3.b) = \\
&  V^{-1}.(a-(\lambda_3+1).b).V^2.(1 + \alpha.b^{p_2})^{-1}.(a - (\lambda_1-2).b).V^{-1} 
\end{align*}
where \ $V = 1 + \beta.b^{p_2}$ \ and  \ $\beta : = (1 + p_2/p_1).\alpha$ \ we see that \ $E$ \ contains a normal rank \ $2$ \ sub-theme which has  fundamental invariants equal to \ $(\lambda_2+1, \lambda_3+1)$ \ and parameter \ $\beta$. $\hfill \square$

\parag{proof of proposition \ref{Uniqueness}} As the result is obvious for \ $k = 1$, we may assume \ $k \geq 2$ \ and the result proved in rank \ $\leq k-1$. Let \ $F_j, j \in [1,k]$ \ and \ $G_j, j \in [1,k]$ \ two J-H. principal sequences for \ $E$. As the sequences \ $\lambda_j + j$ \ and \ $\mu_j+j$ \ co{\"i}ncide up to the order and are both non decreasing, they  co{\"i}ncide. Now let \ $j_0$ \ be the first integer in \ $[1,k]$ \ such that \ $F_{j_0} \not= G_{j_0}$. If \ $j_0 \geq 2$ \ applying the induction hypothesis to \ $E\big/F_{j_0-1}$ \ gives \ $F_{j_0}\big/F_{j_0-1} = G_{j_0}\big/F_{j_0-1}$ \ and so \ $F_{j_0} = G_{j_0}$. \\
So we may assume that \ $j_0 = 1$. Let \ $H$ \ be the normalization of \ $F_1 + G_1$. As \ $F_1$ \ and \ $G_1$ \ are normal rank 1 and distinct, then \ $H$ \ is a rank 2 normal submodule. It is a \ $[\lambda]-$primitive fresco of rank 2 with two normal rank 1 sub-modules which are isomorphic as \ $\lambda_1 = \mu_1$. Moreover the principal J-H. sequence of \ $H$ \ begins by a normal submodule isomorphic to \ $E_{\lambda_1}$. So the previous lemma  implies \ $F_1 = G_1$. So for any \ $j \in [1,k]$ \ we have \ $F_j = G_j$. $\hfill \blacksquare$\\

\begin{defn}\label{Inv. fond.}
Let \ $E$ \ be a \ $[\lambda]-$primitive fresco and consider its principal J-H. sequence \ $F_j, j \in [1.k]$. Put \ $F_j\big/F_{j-1} \simeq E_{\lambda_j}$ \ for \ $j \in [1,k]$ \ (with \ $F_0 = \{0\}$). We shall call {\bf fundamental invariants} of \ $E$ \ the ordered k-tuple\ $(\lambda_1, \dots, \lambda_k) $.
\end{defn}

Of course, if we have any J-H. sequence for a \ $[\lambda]-$primitive fresco  \ $E$ \ with rank \ $1$ \ quotients associated to the rational numbers \ $\mu_1, \dots, \mu_k$, it is easy to recover the fundamental invariants of \ $E$ \ because the numbers \ $\mu_j+j, j \in [1,k]$ \ are the same than the numbers \ $\lambda_j+j$ \ up to a permutation. But as the sequence \ $\lambda_j+j$ \ is non decreasing with \ $j$, it is enough to put the \ $\mu_j+j$ \ in the non decreasing order with \ $j$ \  to conclude.

\bigskip

Let \ $\lambda_1, \dots, \lambda_k$ \ be the fundamental invariants of a \ $[\lambda]-$primitive fresco \ $E_0$, and not \ $\mathcal{F}(\lambda_1, \dots, \lambda_k)$ \ the set of isomorphism classes of frescos with these fundamental invariants. The uniqueness of the principal J-H. sequence of a \ $[\lambda]-$primitive  fresco allows to define for each \ $(i,j) \ 1 \leq i < j \leq k$ \ a map
$$ q_{i,j} : \mathcal{F}(\lambda_1, \dots, \lambda_k) \to \mathcal{F}(\lambda_i, \dots, \lambda_j)$$
defined by \ $q_{i,j}([E]) = [F_j\big/F_{i-1}] $ \ where \ $(F_h)_{h\in [0,k]}$ \ is the principal J-H. sequence of \ $E$. This makes sens because any isomorphism \ $\varphi : E^1 \to E^2$ \ between two \ $[\lambda]-$primitive frescos induces isomorphisms between each term of  the corresponding principal J-H. sequences.\\
For instance classification of rank \ $2$ \ $[\lambda]-$primitive frescos gives (see example before proposition \ref{J-H exist}) for any rational number \ $\lambda_1 > 1$ \ and \ $p_1 \in \mathbb{N}$ :
\begin{align*}
&{\rm for} \quad  p_1 = 0 \quad  \mathcal{F}(\lambda_1, \lambda_1-1) = \{pt\}  \\
& {\rm for} \quad p_1 \geq 1 \quad \alpha : \mathcal{F}(\lambda_1, \lambda_1 + p_1 - 1) \overset{\simeq}{\longrightarrow}  \C 
\end{align*}
and \ $\{pt\}$ \ is given by the isomorphism class of \ $\A\big/\A.(a - \lambda_1.b).(a - (\lambda_1-1).b) $ \ and the isomorphism class associated to \ $\alpha^{-1}(x)$ \ for \ $x  \in \C$ \ in the case \ $p_1 \geq 1$ \ is given by \ $\A\big/\A.(a - \lambda_1.b).(1 + x.b^{p_1})^{-1}.(a - (\lambda_1 + p_1 -1).b)$. \\

\section{Semi-simplicity.}

\subsection{Semi-simple regular  (a,b)-modules.}

\begin{defn}\label{ss 0}
Let \ $E$ \ be a regular (a,b)-module. We say that \ $E$ \ is {\bf semi-simple} if it is a sub-module of a finite direct sum of rank 1 regular (a,b)-modules.
\end{defn}

Note that if \ $E$ \ is a sub-module of a regular (a,b)-module it is necessary regular. As a direct sum of regular (a,b)-modules if regular, our assumption that \ $E$ \ is regular is superfluous.\\
It is clear from this definition that a sub-module of a semi-simple (a,b)-module is semi-simple and that a (finite) direct sum of semi-simple (a,b)-modules is again semi-simple.

\parag{Remark} A rather easy consequence of the classification of rank \ $2$ \ (a,b)-modules given in [B.93]  is that the rank \ $2$ \ simple pole (a,b)-modules defined in the \ $\C[[b]]-$basis \ $x,y$ \  by the relations :
$$(a - (\lambda+p).b).x = b^{p+1}.y \quad {\rm and} \quad (a - \lambda).y = 0$$
for any \ $\lambda \in \C$ \ and \ any \ $p \in \mathbb{N}$ \ are not semi-simple. We leave the verification of this point to the reader.

Let us begin by a characterization of the semi-simple (a,b)-modules which have a simple pole. First we shall prove that a quotient of a semi-simple (a,b)-module is semi-simple. This will be deduced from the following lemma.

\begin{lemma}\label{quot. ss 1}
Let \ $E$ \ be a (a,b)-module which is direct sum of regular rank \ $1$ \ (a,b)-modules, and let \ $F \subset E$ \ be a rank  $1$ \ normal sub-module. Then \ $F$ \ is a direct factor of \ $E$.
\end{lemma}

\begin{cor}\label{quot. ss 2}
If \ $E$ \ is a semi-simple regular (a,b)-module and \ $F$ \ a normal sub-module of \ $E$, the quotient \ $E\big/F$ \ is a (regular) semi-simple (a,b)-module.
\end{cor}

\parag{Proof of the lemma} Let  \ $E = \oplus_{j=1}^k \ E_{\lambda_j}$ \ and assume that \ $F \simeq E_{\mu}$. Let \ $e_j$ \ be a standard generator of \ $E_{\lambda_j}$ \ and \ $e$ \ a standard generator of \ $E_{\mu}$. Write
$$ e = \sum_{j=1}^k \ S_j(b).e_j $$
and compute \ $(a - \mu.b).e = 0$ \ using the fact that \ $e_j, j \in [1,k]$ \ is a \ $\C[[b]]-$basis of \ $E$ \ and the relations \ $(a - \lambda_j.b).e_j = 0$ \ for each \ $j$. We obtain for each \ $j \in [1,k]$ \ the relation
$$ b.S'_j - (\mu-\lambda_j).S_j = 0 .$$
So, if \ $\mu - \lambda_j$ \ is not in \ $\mathbb{N}$, we have \ $S_j = 0$. When \ $\mu = \lambda_j + p_j$ \ with \ $p_j \in \mathbb{N}$ \ we obtain \ $S_j(b) = \rho_j.b^{p_j}$ \ for some \ $\rho_j \in \C$. As we assume that \ $e$ \ is not in \ $b.E$, there exists at least one \ $j_0 \in [1,k]$ \ such that \ $p_{j_0} = 0$ \ and \ $\rho_{j_0} \not= 0 $. Then it is clear that we have
$$ E = F \oplus \big(\oplus_{j\not= j_0} E_{\lambda_j}\big) $$
concluding the proof. $\hfill \blacksquare$

\parag{Proof of the corollary} We argue by induction on the rank of \ $F$. In the rank \ $1$ \ case, we have \ $F \subset E \subset \mathcal{E} : = \oplus_{j=1}^k \ E_{\lambda_j}$. Let \ $\tilde{F}$ \ the normalization of \ $F$ \ in \ $\oplus_{j=1}^k \ E_{\lambda_j}$. Then the lemma shows that there exists a \ $j_0 \in [1,k]$ \ such that 
 $$\mathcal{E} = \tilde{F} \oplus \oplus_{j\not= j_0} E_{\lambda_j}.$$
 Then, as \ $\tilde{F} \cap E = F$, the quotient map \ $\mathcal{E} \to \mathcal{E}\big/\tilde{F} \simeq \oplus_{j\not= j_0} E_{\lambda_j}$ \ induces an injection of \ $E\big/F$ \ in a direct sum of regular rank \ $1$ \ (a,b)-modules. So \ $E\big/F$ \ is semi-simple.\\
 Assume now that the result is proved for \ $F$ \ with rank \ $\leq d-1$ \ and assume that \ $F$ \ has rank \ $d$. Then using a rank \ $1$ \ normal sub-module \ $G$ \ in \ $F$, we obtain that \ $F\big/G $ \ is a normal rank \ $d-1$ \ sub-module of \ $E\big/G$. Using the rank \ $1$ \ case we know that \ $E\big/G$ \ is semi-simple, and the induction hypothesis gives that 
  $$E\big/F = (E\big/G)\Big/(F\big/G)$$
   is semi-simple. $\hfill \blacksquare$\\
   

\begin{prop}\label{simple pole ss}
Let \ $E$ \ be a simple pole semi-simple (a,b)-module. Then \ $E$ \ is  a direct sum of regular rank \ $1$ \ (a,b)-modules.
\end{prop}


\parag{proof}  We shall prove the proposition by induction on the rank \ $k$ \ of \ $E$. As the rank \ $1$ \ case is obvious, assume the proposition proved for \ $k-1$ \ and  that \ $E$ \ is a simple pole rank \ $k \geq 2$ \ semi-simple (a,b)-module. From the existence of J-H. sequence, we may find an exact sequence
$$ 0  \to F \to E \to  E_{\lambda} \to  0 $$
where \ $F$ \ has rank \ $k-1$. By definition \ $F$ \ is semi-simple, but it also has a simple pole because \ $a.F \subset a.E \subset b.E $ \ implies that \ $a.F \subset F \cap b.E = b.F$ \ as \ $F$ \ is normal in \ $E$.
So by the induction hypothesis \ $F$ \ is a direct sum of  regular rank \ $1$ \ (a,b)-modules.\\
Let \ $e \in E$ \ such that its image in \ $E_{\lambda}$ \ is \ $e_{\lambda}$ \ a standard generator of \ $E_{\lambda}$ \ satisfying \ $a, e_{\lambda} = \lambda.b.e_{\lambda}$. Then we have \ $(a - \lambda.b).e \in F$. We shall first look at the case \ $k = 2$. So \ $F$ \ is a rank \ $1$ \ and  we have \ $F \simeq E_{\mu}$ \ for some \ $\mu \in \C$. Let \ $e_{\mu}$ \ be a standard generator in \ $F$ \ and put
$$ (a - \lambda.b).e = S(b).e_{\mu}.$$
 Our simple pole assumption on \ $E$ \ implies \ $S(0) = 0$ \ and we may write \\
  $S(b) = b.\tilde{S}(b)$. We look for \ $T \in \C[[b]]$ \ such that \ $\varepsilon : = e + T(b).e_{\mu}$ \ satisfies
$$ (a - \lambda.b).\varepsilon = 0.$$
If such a \ $T$ \ exists, then we would have \ $E = E_{\mu}\oplus E_{\lambda}$ \ where \ $E_{\lambda} $ \ is the submodule generated by \ $\varepsilon$, because it is clear that we have \ $E = F \oplus \C[[b]].e$ \ as a \ $\C[[b]]-$module. \\
To find \ $T$ \ we have to solve in \ $\C[[b]]$ \  the differential equation
$$ b.\tilde{S}(b) + (\mu - \lambda).b.T(b) + b^2.T'(b) = 0 .$$
If \ $\lambda- \mu$ \ is not a non negative integer, such a \ $T$ \ exists and is unique. But when \ $\lambda = \mu + p$ \ with \ $p \in \mathbb{N}$, the solution exists if and only if the coefficient of \ $b^p$ \ in \ $\tilde{S}$ \ is zero. If it is not the case, define \ $\tilde{T}$ \ as the solution of the differential equation
$$ \tilde{S}(b) -  \alpha.b^p + b.\tilde{T}'(b) - p.\tilde{T}(b) = 0 $$
where \ $\alpha$ \ is the coefficient of \ $b^p$ \ in \ $\tilde{S}$ \ and where we choose \ $\tilde{T}$ \ by asking that it has no \ $b^p$ \ term. Then \ $\varepsilon_1 : = e + \tilde{T}(b).e_{\lambda-p}$ \ satisfies
$$ (a - (\lambda - p).b).\varepsilon_1 = \alpha.b^{p+1}.e_{\lambda -p} .$$
Then, after changing \ $\lambda-p$ \ in \ $\lambda$ \ and  \ $e_{\lambda-p}$ \ in \ $\alpha.e_{\lambda-p}$,  we recognize one of the  rank \ $2$ \ modules  which appears in the previous remark and which is not semi-simple. So we have a contradiction. This concludes the rank \ $2$ \ case.

\bigskip

Now consider the case \ $k \geq  3$ \ and using the induction hypothesis write
$$ F \simeq \oplus_{j=1}^{k-1} \ E_{\mu_j} $$
and denote \ $e_j$ \ a standard generator for \ $E_{\mu_j}$. Write
$$ (a - \lambda.b).e = \sum_{j=1}^{k-1} \ S_j(b).e_j  .$$
The simple pole assumption again gives \ $S_j(0) = 0 $ \ for each \ $j$; now we 
look for \ $T_1, \dots ,T_{k-1}$ \ in \ $\C[[b]]$ \ such that, defining \ $\varepsilon : = e +  \sum_{j=1}^{k-1} \ T_j(b).e_j $, we have
$$ (a - \lambda.b).\varepsilon =  0  .$$
For each \ $j$ \ this leads to the differential equation
$$b.\tilde{S}_j(b) + (\mu_j - \lambda).b.T_j(b) + b^2.T_j'(b) = 0 $$
where we put \ $S_j(b) = b.\tilde{S}_j(b)$. We are back, for each given \ $j \in [1,k-1]$, to the previous problem in the rank \ $2$ \ case. So if for each \ $j \in [1,k-1]$ \ we have a solution \ $T_j \in \C[[b]]$ \  it is easy to conclude that
$$ E \subset F \oplus \C[[b]].\varepsilon \simeq F \oplus E_{\lambda}.$$
If for some \ $j_0$ \ there is no solution, the coefficient of \ $b^{j_0}$ \ in \ $\tilde{S}_{j_0}$ \ does not vanish  and then the image of \ $E$ \ in \ $E\big/\oplus_{j \not= j_0} E_{\mu_j} $ \ is a rank \ $2$ \ not semi-simple (a,b)-module. This contradicts the corollary \ref{quot. ss 2} and concludes the proof. $\hfill \blacksquare$\\

\parag{Remark} Let \ $E$ \ be a semi-simple  rank \ $k$ \ (a,b)-module and let \ $E^{\sharp}$ \ its saturation by \ $b^{-1}.a$, then \ $E^{\sharp}$ \ is semi-simple because there exists \ $N \in \mathbb{N}$ \  with an inclusion \  $E^{\sharp} \subset b^{-N}.E $. Also \ $E^{\flat}$ \ its maximal simple pole sub-module is semi-simple. Then we have 
$$ E^{\flat} \simeq \oplus_{j=1}^k \ E_{\lambda_j+d_j} \subset E \subset E^{\sharp} \simeq \oplus_{j=1}^k \ E_{\lambda_j} $$
where \ $d_1, \dots, d_k$ \ are non negative integers.\\

\begin{cor}\label{dual ss ; quot. ss}
Let \ $E$ \ be a semi-simple (a,b)-module. Then its dual \ $E^*$ \ is semi-simple. 
\end{cor}

\parag{proof} For a regular (a,b)-module \ $E$ \ we have \ $(E^{\flat})^* = (E^*)^{\sharp}$. But the dual of \ $E_{\lambda}$ \ is \ $E_{-\lambda}$ \ and the duality commutes with direct sums. So we conclude that \ $E^* \subset (E^{\flat})^*$ \ and so \ $E^*$ is a semi-simple. $\hfill \blacksquare$

\parag{Remark} The tensor product of two semi-simple (a,b)-modules is semi-simple as a consequence of the fact that \ $E_{\lambda} \otimes E_{\mu} = E_{\lambda+\mu}$.

\subsection{The semi-simple filtration.}

\bigskip

\begin{defn}\label{element ss 1}
Let \ $E$ \ be a regular (a,b)-module and \ $x$ \ an element in \ $E$. We shall say that \ $x$ \ is {\bf semi-simple} if \ $\A.x$ \ is a semi-simple (a,b)-module.
\end{defn}

It is clear that any element in a semi-simple  (a,b)-module is semi-simple. The next lemma shows that the converse is true.

\begin{lemma}\label{element ss 2}
Let \ $E$ \ be a regular (a,b)-module such that any \ $x \in E$ \ is semi-simple. Then \ $E$ \ is semi-simple.
\end{lemma}

\parag{proof} Let \ $e_1, \dots, e_k$ \ be a \ $\C[[b]]-$basis of \ $E$. Then each \ $\A.e_j$ \ is semi-simple, and the map \ $\oplus_{j=1}^k \A.e_j \to E $ \ is surjective. So \ $E$ \ is semi-simple thanks to the corollary \ref{quot. ss 2} and the comment following the definition \ref{ss 0}. $\hfill \blacksquare$\\

\begin{lemma}\label{ss part 0}
Let \ $E$ \ be a regular (a,b)-module. The subset \ $S_1(E)$ \ of semi-simple elements in \ $E$ \ is a normal submodule in \ $E$.
\end{lemma}

\parag{proof} As a direct sum and quotient of semi-simple (a,b)-modules are semi-simple, it is clear that for  \ $x$ \ and \ $y$ \ semi-simple the sum \ $\A.x + \A.y$ \ is semi-simple. So \ $x + y$ \ is semi-simple. This implies that \ $S_1(E)$ \ is a submodule of \ $E$. If \ $b.x$ \ is in \ $S_1(E)$, then \ $\A.b.x$ \ is semi-simple. Then \ $\A.x \simeq \A.b.x \otimes E_{-1}$ \ is also semi-simple, and so \ $S_1(E)$ \ is normal in \ $E$. $\hfill \blacksquare$

\begin{defn}\label{ss part 1}
Let \ $E$ \ be a regular (a,b)-module. Then the submodule \ $S_1(E)$ \ of semi-simple elements in \ $E$ \ will be called the {\bf semi-simple part} of \ $E$.\\
Defining \ $S_j(E)$ \ as the pull-back on \ $E$ \ of the semi-simple part of \ $E\big/ S_{j-1}(E)$ \ for \ $j \geq 1$ \ with the initial condition \ $S_0(E) = \{0\}$, we define a sequence of normal submodules in \ $E$ \ such that \ $S_j(E)\big/S_{j-1}(E) = S_1(E\big/S_{j-1}(E)$. We shall call it the {\bf \em semi-simple filtration} of \ $E$. The smallest integer \ $d$ \ such we have \ $S_d(E) = E$ \ will be called the {\bf semi-simple depth} of \ $E$ \ and we shall denote it \ $d(E)$.
\end{defn}

\parag{Remarks}
\begin{enumerate}[i)]
\item As \ $S_1(E)$ \ is the maximal semi-simple sub-module of \ $E$ \ it contains any rank \ $1$ \ sub-module of \ $E$. So \ $S_1(E) = \{0\}$ \ happens if and only if \ $E = \{0\}$.
\item Then the ss-filtration of \ $E$ \ is strictly increasing for \ $0 \leq j \leq d(E)$. 
\item It is easy to see that for any submodule \ $F$ \ in \ $E$ \ we have \ $S_j(F) = S_j(E) \cap F$. So \ $S_j(E)$ \ is the subset of \ $x \in E$ \ such that \ $ d(x) : = d(\A.x) \leq j$ \ and \  $d(E)$ \ is the supremum of \ $d(x)$ \ for \ $x \in E$.  $\hfill \square$
\end{enumerate}

\begin{lemma}\label{primitive ss}
Let \ $0 \to F \to E \to G \to 0$ \ be an exact sequence of regular (a,b)-modules such that \ $F$ \ and \ $G$ \ are semi-simple and such that for any roots \ $\lambda, \mu$ \ of the Bernstein polynomial of \ $F$ \ and \ $G$ \ respectively we have \ $\lambda - \mu \not\in \mathbb{Z}$. Then \ $E$ \ is semi-simple.
\end{lemma}

\parag{proof} By an easy induction we may assume that \ $F$ \ is rank \ $1$ \ and also that \ $G$ \ is rank \ $1$. Then the conclusion follows from the classification of the rank \ $2$ \ regular (a,b)-module. $\hfill \blacksquare$

 \parag{Application} A regular (a,b)-module \ $E$ \ is semi-simple if and only if for each \ $\lambda \in \C$ \ its primitive part \ $E[\lambda]$ \ is semi-simple.\\
   
The following easy facts are  left as  exercices for the reader.

\begin{enumerate}
\item  Let \ $0 \to F \to E \to G \to 0$ \ be a short exact sequence of regular (a,b)-modules. Then we have the inequalities
  $$ \sup \{ d(F), d(G)\} \leq d(E) \leq d(F) + d(G) .$$
  \item Let \ $E \subset \Xi_{\lambda}^{(N)} \otimes V$ \ be a sub-module where \ $V$ \ is  a finite dimensional complex space. Then for each \ $j \in [1, N+1]$ \ we have
  $$ S_j(E) = E \cap \left( \Xi_{\lambda}^{(j-1)} \otimes V \right) .$$
\item  In the same situation as above,  the minimal \ $N$ \ for such an embedding of a geometric \ $[\lambda]-$primitive \ $E$ \ is equal to \ $d(E)-1$.
\item In the same situation as above, the equality  \ $S_1(E) = E \cap \Xi_{\lambda}^{(0)} \otimes V$ \ implies \ $rank(S_1(E)) \leq \dim_{\C} V $. It is easy to see that any \ $[\lambda]-$primitive semi-simple geometric (a,b)-module  \ $F$ \ may be embedded in \ $\Xi_{\lambda}^{(0)} \otimes V$ \ for some \ $V$ \ of dimension \ $rank(F)$. But it is also easy to prove that when \ $F = S_1(E)$ \ with \ $E$ \ as above, such an embedding may be extended to an embedding of \ $E$ \ into \ $\Xi_{\lambda}^{(d(E)-1)} \otimes V$. $\hfill \square$
\end{enumerate}

\subsection{Semi-simple frescos.}

We begin by a simple remark : A  \ $[\lambda]-$primitive theme is semi-simple if and only if it has rank \ $\leq 1$. This is an easy consequence of the fact that the saturation of a rank \ $2$ \ $[\lambda]-$primitive theme is one of the rank \ $2$ \ (a,b)-modules considered in the remark following the definition \ref{ss 0} which are not semi-simple

\begin{lemma}\label{semi-simple}
A  geometric (a,b)-module \ $E$ \ is {\bf semi-simple} if and only if any  quotient of \ $E$ \ which is a \ $[\lambda]-$primitive theme for some \ $[\lambda] \in \mathbb{Q}\big/\mathbb{Z}$ \   is of rank \ $\leq 1$.
\end{lemma}

\parag{proof} The condition is clearly necessary as a quotient of a semi-simple (a,b)-module is semi-simple (see corollary \ref{quot. ss 2}), thanks to the remark above. Using the application of the lemma 
\ref{primitive ss}, it is enough to consider the case of a \ $[\lambda]-$primitive fresco to prove that the condition is sufficient. Let \ $\varphi : E \to \Xi_{\lambda}^{(N)} \otimes V$ \ be an embedding of \ $E$ \ which exists thanks to the embedding theorem 4.2.1. of [B.09], we obtain that each component of this map in a basis \ $v_1, \dots v_p$ \ of \ $V$ \ has rank at most \ $1$ \ as its image is a \ $[\lambda]-$primitive theme. Then each of these images is isomorphic to some \ $E_{\lambda+q}$ \ for some integer \ $q$. So we have in fact an embedding of \ $E$ \  in a direct sum of \ $E_{\lambda+ q_i}$ \ and \ $E$ \ is semi-simple. $\hfill \blacksquare$

\parag{Remark} The preceeding proof shows that a fresco is a non zero  \ $[\lambda]-$primitive theme for some \ $[\lambda] \in \mathbb{Q}\big/\mathbb{Z}$ \ if and only if \ $S_1(E)$ \ has rank  $1$. $\hfill \square$

\begin{lemma}\label{all J-H}
Let \ $E$ \ be a semi-simple fresco  with rank \ $k$ \ and let \ $\lambda_1, \dots, \lambda_k$ \ be the numbers associated to a J-H. sequence of \ $E$. Let \ $\mu_1, \dots, \mu_k$ \ be a twisted permutation\footnote{This means that the sequence \ $\mu_j+j, j \in [1,k]$ \ is a permutation (in the usual sens) of \ $\lambda_j+j , j\in [1,k]$.} of \ $\lambda_1, \dots, \lambda_k$. Then there exists a J-H. sequence for \ $E$ \ with quotients corresponding to \ $\mu_1, \dots, \mu_k$.
\end{lemma}

\parag{Proof} As the symetric group \ $\mathfrak{S}_k$ \  is generated by the transpositions \ $t_{j,j+1}$ \ for \ $j \in [1,k-1]$, it is enough to show that, if \ $E$ \ has a J-H. sequence with quotients given by the numbers \ $\lambda_1, \dots, \lambda_k$, then there exists a J-H. sequence for \ $E$ \ with quotients \ $\lambda_1, \dots, \lambda_{j-1}, \lambda_{j+1}+1,\lambda_j-1, \lambda_{j+2}, \dots, \lambda_k$ \ for \ $j \in [1,k-1]$. But \ $G : = F_{j+1}\big/F_{j-1}$ \ is a rank 2  sub-quotient of \ $E$ \ with an exact sequence
$$ 0 \to E_{\lambda_j} \to G \to E_{\lambda_{j+1}} \to 0 .$$
As \ $G$ \ is  a rank 2 semi-simple fresco, it admits also an exact sequence
$$ 0 \to G_1 \to G \to G\big/G_1 \to 0$$
with \ $G_1 \simeq E_{\lambda_{j+1}+1}$ \ and \ $G\big/G_1 \simeq E_{\lambda_j-1}$. Let \ $q : F_{j+1} \to G$ \ be the quotient map. Now the J-H. sequence for \ $E$ \ given by
$$ F_1, \dots, F_{j-1}, q^{-1}(G_1), F_{j+1}, \dots ,F_k = E$$
 satisfies our requirement. $\hfill \blacksquare$
 
  \begin{prop}\label{crit. ss}
 Let \ $E$ \ be a  \ $[\lambda]-$primitive fresco. A necessary and sufficient condition in order that \ $E$ \ is semi-simple is that it admits a J-H. sequence with quotient corresponding to \ $\mu_1, \dots, \mu_k$ \ such that the sequence \ $\mu_j+j$ \ is strictly decreasing.
 \end{prop}
 
 \parag{Remarks}
 \begin{enumerate}
 \item As a fresco is semi-simple if and only if for each \ $[\lambda]$ \ its  \ $[\lambda]-$primitive  part is semi-simple, this proposition gives also a criterium to semi-simplicity for any fresco.
  \item This criterium is a very efficient tool to produce easily examples of semi-simple frescos.
  \end{enumerate}
 
 \parag{Proof} Remark first that if we have, for a  \ $[\lambda]-$primitive fresco \ $E$, a J-H. sequence \ $F_j, j \in [1,k]$ \  such that \ $\lambda_j+j = \lambda_{j+1}+j+1$ \ for some \ $j \in [1,k-1]$, then \ $F_{j+1}\big/F_{j-1}$ \ is a sub-quotient of \ $E$ \ which is a  \ $[\lambda]-$primitive theme of rank  2. So \ $E$ \ is not semi-simple. As a consequence, when a \ $[\lambda]-$primitive fresco \ $E$ \ is semi-simple the principal J-H. sequence corresponds to a strictly increasing sequence \ $\lambda_j+j$. Now, thanks to the previous lemma we may find a J-H. sequence for \ $E$ \  corresponding to the strictly decreasing order for the sequence \ $\lambda_j+j$. \\

 No let us prove the converse. We shall use the following lemma.
 
 \begin{lemma}\label{utile}
 Let \ $F$ \ be a rank \ $k$ \ semi-simple  \ $[\lambda]-$primitive fresco and let \ $\lambda_j+j$ \ the strictly increasing sequence corresponding to its principal J-H. sequence. Let \ $\mu \in [\lambda]$ \ such that \ $0 < \mu+k < \lambda_1+1$. Then any fresco \ $E$ \ in an exact sequence
 $$ 0 \to F \to E \to E_{\mu} \to 0 $$
 is semi-simple (and \ $[\lambda]-$primitive).
 \end{lemma}
 
 \parag{Proof} The case \ $k=1$ \ is obvious, so assume that \ $k \geq 2$ \ and that we have a rank 2 quotient \ $\varphi : E \to T$ \ where \ $T$ \ is a \ $[\lambda]-$primitive  theme. Then \ $Ker\, \varphi\, \cap F$ \ is a normal submodule of \ $F$ \ of rank \ $k-2$ \ or \ $k-3$ \ (for \ $k \geq 3$). If \ $Ker\,\varphi \,\cap F$ \ is of rank \ $k-3$, the rank of \ $F\big/(Ker\, \varphi\, \cap F)$ \ is \ $2$ \  and it injects in \ $T$ \ via \ $\varphi$. So \ $F\big/(Ker\, \varphi \cap F)$ \ is a rank 2 \ $[\lambda]-$primitive theme. As it is semi-simple, because \ $F$ \ is semi-simple, we get a contradiction.\\
 So the rank of \ $F\big/(Ker\, \varphi\, \cap F)$ \ is \ $1$ \ and we have an exact sequence
 $$ 0 \to F\big/(Ker\, \varphi\, \cap F) \to T \to E\big/F \to 0 .$$
 Put \ $F\big/(Ker\, \varphi\, \cap F) \simeq F_1(T)  \simeq E_{\nu}$. Because \ $T$ \ is a \ $[\lambda]-$primitive  theme, we have the inequality \ $\nu+1 \leq \mu+2$. Looking at a J-H. sequence of \ $E$ \ ending by
 $$ \dots \subset Ker\, \varphi \cap F \subset  \varphi^{-1}(F_1(T)) \subset E $$
 we see that \ $\nu+k-1$ \ is in the set \ $\{ \lambda_j+j, j \in [1,k]\}$ \  and as \ $\lambda_1+1$ \ is the infimum of this set we obtain  \ $\lambda_1+1 \leq \mu+k $ \ contradicting our assumption. $\hfill \blacksquare$
 
 \parag{End of proof of the proposition \ref{crit. ss}} Now we shall prove by induction on the rank of a \ $[\lambda]-$primitive fresco \ $E$ \ that if it admits a J-H. sequence corresponding to a strictly decreasing sequence \ $\mu_j +j$, it is semi-simple. As the result is obvious in rank 1, we may assume \ $k \geq 1$ \ and the result proved for \ $k$. So let \ $E$ \ be a fresco of rank \ $k+1$ \ and let \ $F_j, j\in [1,k+1]$ \ a J-H. sequence for \ $E$ \ corresponding to the strictly decreasing sequence \ $\mu_j+j, j \in [1,k+1]$. Put \ $F_{j}\big/F_{j-1} \simeq E_{\mu_j}$ \ for all \ $j \in [1,k+1]$, define \ $F : = F_k$ \ and \ $\mu : = \mu_{k+1}$; then the induction hypothesis gives that \ $F$ \ is semi-simple and we apply the previous lemma to conclude. $\hfill \blacksquare$\\

    The following interesting corollary is an obvious consequence of the previous proposition.
  
 \begin{cor}\label{ss base 3}
  Let \ $E$ \ be a   fresco and let \ $\lambda_1, \dots, \lambda_k$ \ be the numbers associated to any J-H. sequence of \ $E$. Let \ $\mu_1, \dots, \mu_d$ \ be the numbers associated to any  J-H. sequence of \ $S_1(E)$. Then, for \ $j \in [1,k]$,  there exists a rank 1 normal submodule of \ $E$ \ isomorphic to \ $E_{\lambda_j+j-1}$ \ if and only if there exists \ $i \in [1,d]$ \ such that we have  \ $\lambda_j+j-1 = \mu_i+i-1$.
  \end{cor}

 Of course, this gives the list of all isomorphy classes of  rank 1 normal submodules  contained in \ $E$. So, using shifted duality, we get also the list of all isomorphy classes of  rank 1 quotients of \ $E$. It is interesting to note that this is also the list of the possible initial exponents for maximal logarithmic terms which appear in the asymtotic expansion of a given relative de Rham cohomology class after integration on any vanishing cycle in the spectral subspace of the monodromy associated to the eigenvalue \ $exp(2i\pi.\lambda)$.\\

\section{A criterion for semi-simplicity.}

\subsection{Polynomial dependance.}

All \ $\C-$algebras have a unit.\\

When we consider a sequence of algebraically independent variables \ $\rho : = (\rho_i)_{i \in \mathbb{N}}$ \ we shall denote by \ $\C[\rho]$ \ the \ $\C-$algebra generated by these variables. Then \ $\C[\rho][[b]]$ \ will be the commutative \ $\C-$algebra of formal power series
$$ \sum_{\nu = 0}^{\infty} \ P_{\nu}(\rho).b^{\nu} $$
where \ $P_{\nu}(\rho)$ \ is an element in \ $\C[\rho]$ \ so a polynomial in \ $\rho_0, \dots, \rho_{N(\nu)} $ \ where \ $N(\nu)$ \ is an integer depending on \ $\nu$. So each coefficient in the formal power serie in \ $b$ \ depends only on a finite number of the variables \ $\rho_i$.

\begin{defn}\label{polyn. depend.}
Let \ $E$ \ be a (a,b)-module and let \ $e(\rho)$ \ be a family of elements in  \ $E$ \ depending on a family of variables \ $(\rho_i)_{i\in \mathbb{N}}$. We say that \ $e(\rho)$ \ {\bf\em depends polynomially on \ $\rho$} \ if there exists a fixed \ $\C[[b]]-$basis \ $e_1, \dots, e_k$ \ of \ $E$ \ such that
$$ e(\rho) = \sum_{j=1}^k \ S_j(\rho).e_j $$
where \ $S_j$ \ is for each \ $j \in [1,k]$ \ an element in the algebra  \ $\C[\rho][[b]]$.
\end{defn}

\parag{Remarks} \begin{enumerate}
\item It is important to remark that when \ $e(\rho)$ \ depends polynomially of \ $\rho$, then \ $a.e(\rho)$ \ also. Then for any \ $u \in \A$, again \ $u.e(\rho)$ \ depends polynomially on \ $\rho$. 
\item It is easy to see that we obtain an equivalent condition on the family \ $e(\rho)$ \ by asking that the coefficient of \ $e(\rho)$ \ are in \ $\C[\rho][[b]]$ \ in a \ $\C[[b]]-$basis  of \ $E$ \  whose elements depend polynomially of \ $\rho$.
\item The invertible elements in the algebra \ $\C[\rho][[b]]$ \ are exactly those elements with a constant term in \ $b$ \ invertible in the algebra \ $\C[\rho]$.  As we assume that the variables \ $(\rho)_{i\in \mathbb{N}}$ \ are algebraically independent, the invertible elements are  those  with a constant term in \ $b$ \ in \ $\C^*$.
 $\hfill \square$\\
\end{enumerate}

\begin{prop}\label{gen. dep.}
Let \ $E$ \ be a rank \ $k$ \ $[\lambda]-$primitive fresco and let \ $\lambda_1, \dots, \lambda_k$ \ its fundamental invariants. Let \ $e(\rho)$ \ be a family of generators of \ $E$ \ depending polynomially on a family of algebraically independant variables \ $(\rho_i)_{i \in \mathbb{N}}$. Then there exists \ $S_1, \dots, S_k$ \ in \ $\C[\rho][[b]]$ \ such that \ $S_j(\rho)[0] \equiv 1$ \ for each \ $j \in [1,k]$ \ and such that the annihilator of \ $e(\rho)$ \ in \ $E$ \ is generated by the element of \ $\A$
$$ P(\rho) : = (a - \lambda_1.b).S_1(\rho)^{-1} \dots S_{k-1}(\rho)^{-1}.(a - \lambda_k.b).S_k(\rho)^{-1} .$$
\end{prop}

\parag{Proof} The key result to prove this proposition is the rank \ $1$ \ case. In this case we may consider a standard generator \ $e_1$ \ of \ $E$ \ which is a \ $\C[[b]]-$basis of \ $E$ \ and satisfies
$$ (a - \lambda_1.b).e_1 = 0 .$$
Then, by definition, we may write
$$ e(\rho) = S_1(\rho).e_1 $$
where \ $S_1 $ \ is in \ $\C[\rho][[b]]$ \ is invertible in this algebra, so has a  constant term in \ $\C^*$. Up to normalizing \ $e_1$, we may assume that \ $S_1(\rho)[0] \equiv 1 $ \ and then define 
 $$P(\rho) : = (a - \lambda_1.b).S_1(\rho)^{-1} .$$
It clearly generates the annihilator of \ $e(\rho)$ \ for each \ $\rho$.\\
Assume now that the result is already proved for the rank \ $k-1 \geq 1$. Then consider the family \ $[e(\rho)]$ \ in the quotient \ $E\big/F_{k-1}$ \ where \ $F_{k-1}$ \ is the rank \ $k-1$ \ sub-module of \ $E$ \ in its principal J-H. sequence. Remark first that \ $[e(\rho)]$ \ is a family of generators of
 \ $E\big/F_{k-1}$ \ which depends polynomially on \ $\rho$. This is a trivial consequence of the fact that we may choose a \ $\C[[b]]-$basis \ $e_1, \dots, e_k$ \  in \ $E$ \ such that \ $e_1, \dots, e_{k-1}$ \ is a \ $\C[[b]]$ \ basis of \ $F_{k-1}$ \ and \ $e_k$ \ maps to a standard generator of \ $E\big/F_{k-1}\simeq E_{\lambda_k}$. Then the rank \ $1$ \  case gives \ $S_k \in \C[\rho][[b]]$ \ with \ $S_k(\rho)[0] = 1$ \  and such that \ $(a - \lambda_k.b).S_k(\rho)^{-1}.e(\rho)$ \ is in \ $F_{k-1}$ \ for each \ $\rho$. But then it is a family of generators of \ $F_{k-1}$ \ which depends polynomially on \ $\rho$ \ and the inductive assumption allows to conclude. $\hfill \blacksquare$\\

Fix now the fundamental invariants \ $\lambda_1, \dots, \lambda_k$ \ for a \ $[\lambda]-$primitive fresco. 

\begin{defn}\label{polyn. depend. E}
Consider now a complex valued  function \ $f$ \  defined on a subset \ $\mathcal{F}_0$ \ of the isomorphism classes \ $\mathcal{F}(\lambda_1, \dots, \lambda_k)$ \ of \ $[\lambda]-$primitive frescos with fundamental invariants \ $\lambda_1, \dots, \lambda_k$. We shall say that  \ $f$ \  {\bf depends polynomially on the isomorphism class  \ $[E] \in \mathcal{F}_0 $ \ of the fresco \ $E$} \ if the following condition is satisfied :\\
Let \ $s$ \ be the collection of algebraically independant variables corresponding to the non constant  coefficients of \ $k$ \ \'elements \ $S_1, \dots, S_k$ \ in \ $\C[[b]]$ \ satisfying \ $S_j(0) = 1\quad \forall j \in [1,k]$ \ and consider for each value of \ $s$ \ the rank \ $k$ \ $[\lambda]-$primitive  fresco \ $E(s) : = \A \big/ \A.P(s) $ \ where 
 $$P(s) : =  (a - \lambda_1.b)S_1(s)^{-1} \dots (a - \lambda_k.b).S_k(s)^{-1} $$
 where \ $S_1(s), \dots, S_k(s)$ \ correspond to the given values for \ $s$.\\
 Then there exists a polynomial \ $F \in \C[s]$ \ such that for each value of  \ $s$ \ such that \ $[E(s)]$ \ is in \ $\mathcal{F}_0$, the value of \ $F(s)$ \ is equal to\ $f([E(s)])$.
 \end{defn}

\parag{Example} Let \ $(s_i)_{i \in \mathbb{N}^*}$ \ a family of algebraically independant variables and let  \ $S(s) : = 1 + \sum_{i=1}^{\infty} \ s_i.b^i \in \C[s][[b]]$. Define \ $E(s) : = \A\big/\A.(a - \lambda_1.b).S(s)[b]^{-1}.(a - \lambda_2.b) $ \ where \ $\lambda_1 > 1$ \ is rational and \ $\lambda_2 : = \lambda_1 + p_1 -1$ \ with \ $p_1 \in \mathbb{N}^*$. Define \ $\alpha(s) : = s_{p_1}$. Then the number \ $\alpha(s)$ \ depends only of the isomorphism class of the fresco \ $E(s)$ \ and defines a function on \ $\mathcal{F}(\lambda_1, \lambda_2)$ \ which depends polynomially on \ $[E] \in \mathcal{F}(\lambda_1, \lambda_2)$.  $\hfill \square$

\subsection{The \ $\alpha-$invariant.}

The first proposition will be the induction step in the construction of the \ $\alpha-$invariant.

\begin{prop}\label{Gk}
Fix \ $k \geq 2$. Denote \ $\mathcal{F}_0(\lambda_1, \dots, \lambda_k)$ \ the subset of \ $\mathcal{F}(\lambda_1, \dots, \lambda_k)$ \ of isomorphism class of rank \ $k$ \ $[\lambda]-$primitive fresco \ $E$ \  with invariants \ $\lambda_1, \dots, \lambda_k$ \ such that \ $F_{k-1}$ \ and \ $E\big/F_1$ \ are semi-simple where \ $(F_j)_{j\in [1,k]} $ \ is  the principal J-H. sequence of \ $E$.\\
We assume that, for \ $k \geq 3$ \ we have :
\begin{itemize}
\item  for  rank \ $\leq k-1$ \ frescos the  \ $\alpha$ \ invariant is defined on the corresponding subset  \ $\mathcal{F}_0$;
\item  for \ $[E] \in \mathcal{F}_0$\ the number  \ $\alpha(E)$ \ is zero if and only if \ $E$ \  is semi-simple 
\item and that \ $[E] \mapsto \alpha(E)$ \   depends polynomially of \ $[E] \in \mathcal{F}_0$.
\end{itemize}
Fix \ $E$ \ with \ $[E]$ \ in \ $\mathcal{F}_0(\lambda_1, \dots, \lambda_k)$. Let \ $e$ \ be a generator of \ $E$ \ such that \ $(a - \lambda_{k-1}.b).(a - \lambda_k.b).e$ \ lies in \ $F_{k-2}$. Then the sub-module \ $G_{k-1}$ \ of \ $E$ \ generated by \ $(a - (\lambda_{k-1}-1).b).e$ \ is a normal rank \ $k-1$ \ sub-module of \ $E$ \ which is in \ $\mathcal{F}_0(\lambda_1, \dots, \lambda_{k-2}, \lambda_k+1)$ \  and \ $\alpha(G_{k-1})$ \ is independant of the choice of such a generator \ $e$. Moreover, it defines a polynomial function on \ $\mathcal{F}_0(\lambda_1, \dots, \lambda_k)$.
\end{prop}

\parag{Proof} We shall prove first that if \ $e$ \ is a generator of \ $e$ \ such that \\
 $(a - \lambda_{k-1}.b).(a - \lambda_k.b).e$ \ lies in \ $F_{k-2}$ \ then \ $G_{k-1} : = \A.(a - (\lambda_{k-1}-1).b.e$ \ is a normal rank \ $k-1$ \ sub-module of \ $E$ \ such that \ $[G_{k-1}] $ \ belongs to \ $\mathcal{F}_0(\lambda_1, \dots, \lambda_{k-2}, \lambda_k+1)$.\\
Using the identity in \ $\A$ :
 $$(a - \lambda_{k-1}.b).(a - \lambda_k.b) = (a - (\lambda_k+1).b).(a - (\lambda_{k-1}-1).b) $$
 we see that \ $G_{k-1}$ \ contains \ $F_{k-2}$, that \ $G_{k-1}\big/F_{k-2} \simeq E_{\lambda_k+1}$ \  and \ $G_{k-1}$ \ admits 
 $$ F_1 \subset \dots \subset F_{k-2} \subset G_{k-1} $$
 as principal J-H. sequence. Then the  fundamental invariants for \ $G_{k-1}$ \ are equal to \ $\lambda_1, \dots, \lambda_{k-2},\lambda_k+1$, and the corank \ $1$ \  term \ $F_{k-2}$ \ is semi-simple. As \ $G_{k-1}\big/F_1$ \ is a sub-module of \ $E\big/F_1$ \ it is semi-simple and we have proved that 
 $$ [G_{k-1}] \in \mathcal{F}_0(\lambda_1, \dots, \lambda_{k-2}, \lambda_k+1).$$
 \parag{Claim} If \ $\varepsilon$ \ is another generator of \ $E$ \ such that \ $(a - \lambda_{k-1}.b).(a - \lambda_k.b).\varepsilon $ \ is in \ $F_{k-2}$ \ we have
$$ \varepsilon = \rho.e + \sigma.b^{p_{k-1}-1}.(a - \lambda_k.b).e \quad {\rm modulo} \ F_{k-2} $$
where \ $(\rho, \sigma)$ \ is in \ $\C^*\times \C$.\\
Write \ $\varepsilon = U.e_k + V.e_{k-1} \quad {\rm modulo} \ F_{k-2}$ \ where \ $e_k : = e$ \ and \ $e_{k-1} : = (a - \lambda_k.b).e $ \ and where \ $U, V $ \ are in \ $\C[[b]]$. Now 
$$ (a - \lambda_k.b).\varepsilon = b^2.U'.e_k + U.e_{k-1} + (\lambda_{k-1} - \lambda_k).b.V.e_{k-1} +  b^2.V'.e_{k-1} \quad {\rm modulo} \ F_{k-2} $$
and, as \ $(a - \lambda_{k-1}.b).b^2.U'e_k \in F_{k-1}$ \ implies \ $U' = 0$,  we have  \ $U = \rho \in \C^*$, and we obtain
$$ (a - \lambda_k.b).\varepsilon = \left[ \rho + b^2.V' - (p_{k-1}-1).b.V\right].e_{k-1} \quad {\rm modulo} \ F_{k-2} .$$
If \ $T : = \rho + b^2.V' - (p_{k-1}-1).b.V$ \ we have
$$ (a - \lambda_{k-1}.b).(a - \lambda_k.b).\varepsilon = b^2.T'.e_{k-1}  \quad {\rm modulo} \ F_{k-2} $$
and so  \ $T $ \ is a constant and equal to \ $\rho$. Then \ $V = \sigma.b^{p_{k-1}-1}$ \ for some complex number \ $\sigma$ \ and our claim is proved.

\smallskip

 For \ $ \tau \in \C$ \ define
$$ \varepsilon(\tau) : = e_k + \tau.b^{p_{k-1}-1}.e_{k-1} $$
and let 
$$ G^{\tau}_{k-1} : = \A.(a - (\lambda_{k-1}-1).b).(e_k + \tau.b^{p_{k-1}-1}.e_{k-1}) = \A.(a - (\lambda_{k-1}-1).b).\varepsilon(\tau).$$
Remark also that for any choice of \ $\varepsilon$ \ such that \ $(a -\lambda_{k-1}.b).(a - \lambda_k.b).\varepsilon $ \ is in \ $F_{k-2}$ \ the sub-module genrated by \ $\varepsilon$ \  is equal to \ $G^{\tau}_{k-1}$ \ for some \ $\tau \in \C$.\\
The generator \ $(a - (\lambda_k+1).b).(a - (\lambda_{k-1}-1).b).\varepsilon(\tau) $ \ of \ $F_{k-2}$ \ depends polynomially on \ $\tau$. Then the proposition \ref{gen. dep.} gives \ $Q(\tau)$ \ depending polynomially on \ $\tau$ \ which annihilates for each \ $\tau$ \ this generator. Then \ $Q(\tau).(a - (\lambda_k+1).b) $ \ annihilates the generator \ $(a - (\lambda_{k-1}-1).b).\varepsilon(\tau) $ \ of \ $G^{\tau}_{k-1}$.  Remark now that for any \ $\tau \in \C$ \ the fresco \ $G^{\tau}_{k-1}$ \ has its principal J-H. sequence given by \ $F_1\subset \dots \subset F_{k-2} \subset G^{\tau}_{k-1} $. So \ $G^{\tau}_{k-1}\big/F_1$ \ and \ $F_{k-2}$ \ are  semi-simple, and by our induction hypothesis, the number \ $\alpha(G^{\tau}_{k-1} )$ \ is well defined and \ $\tau \mapsto \alpha(G^{\tau}_{k-1})$ \ is a polynomial in \ $\tau$.\\

The key point will be now to show that either \ $ \alpha(G^{\tau}_{k-1})$ \  is identically zero or it never vanishes. This will prove that this number is independant of \ $\tau$ \ concluding the proof.\\
But if for some \ $\tau \in \C$ \ we have \ $\alpha(G^{\tau}) = 0 $, then \ $G^{\tau}$ \ is semi-simple. Now the sub-module \ $F_{k-1} + G^{\tau}_{k-1} $ \ has  rank \ $k$ \  and  is semi-simple. So its normalization is also semi-simple. But this normalization is \ $E$. So any sub-module of \ $E$ \ is semi-simple and we have \ $\alpha(G^{\tau'}_{k-1}) = 0$ \ for any \ $\tau'$.\\
This also shows that when \ $\alpha(G^{\tau}_{k-1}) \not= 0$ \ for some \ $\tau$ \ then it is not zero for any choice of \ $\tau$. Now a polynomial never vanishing is a constant, and this implies that we may define \ $\alpha(E) : = \alpha(G^{\tau}_{k-1}) $ \ for any value of \ $\tau$.\\
To conclude the proof we have to show that, with this definition, the number \ $\alpha(E)$ \ depends polynomially on \ $E$. For that purpose it is enough, thanks to the previous computation, to produce a polynomial function \ $F \in \C[s]$ \ such that for each \ $s$ \ with \ $E(s) \in \mathcal{F}_0(\lambda_1, \dots, \lambda_k)$ \ we have \ $F(s) = \alpha([E(s)]$.\\
The first remark is that we may assume that \ $S_k = 1$. Then, thanks to the previous computation and the induction on the rank, it is enough to show that we may also reduce to the case where \ $S_{k-1} = 1$ \ because in this case we have 
$$ \alpha(E[s]) = \alpha([G_{k-1}(s)]) $$
where  
$$G_{k-1}(s) = \A\big/\A. (a - \lambda_1.b).S_1(s)^{-1}\dots S_{k-2}(s)^{-1}.(a - (\lambda_k+1).b) .$$
Denote \ $s_{k-1}^i, i \in \mathbb{N}^*$ \ the variables associated to the non constant coefficients of \ $S_{k-1}$.
Looking to a generator of the form
$$ \tilde{e} : = e_k + X.e_{k-1} $$
with \ $X \in \C[s_{k-1}][[b]]$ \ such that \ $(a - \lambda_{k-1}.b).(a - \lambda_k.b).\tilde{e} $ \ is in \ $F_{k-2}$ \ this leads to the equation
\begin{align*}
& b^2.X' - (p_{k-1}-1).b.X = 1 - S_{k-1}(s)
\end{align*}
which has an unique solution \ $X \in \C[s_{k-1}][[b]]$ \ without term in \ $b^{p_{k-1}-1}$ \ as we assume \ $E(s)\big/F_1(s)$ \ semi-simple, hypothesis which implies that \ $S_{k-1}$ \ has no \ $b^{p_{k-1}}$ \ term.
With this generator \ $\tilde{e}$ \ we have \ $(a - \lambda_{k-1}.b).(a - \lambda_k.b).\tilde{e} $ \  which is a generator of \ $F_{k-2}$ \ which depends polynomially on \ $s_{k-1}^i$. Now we may consider \ $F_{k-2}$ \ as a fix fresco (so we fix the coefficients in \ $S_1, \dots, S_{k-2}$ \ and apply the proposition \ref{gen. dep.} to obtain a \ $Q(s) : = (a - \lambda_1.b).S_1(s)^{-1} \dots (a - \lambda_{k-2}.b).S_{k-2}(s)^{-1} $ \ which annihilates \ $(a - \lambda_{k-1}.b).(a - \lambda_k.b).\tilde{e} $ \ in \ $E(s)$ \ and depends polynomially on \ $s$. We may conclude by the induction assumption because we know that we have
$$ \alpha([E(s)]) = \alpha([G_{k-1}(s)]) $$
where \ $G_{k-1}(s) = \A.(a - \lambda_{k-1}-1).b).\tilde{e}(s) $ \ and because we know that the generator \ $(a - \lambda_{k-1}-1).b).\tilde{e}(s) $ \ is annihilated by \ $Q(s).(a - (\lambda_k+1).b)$ \ which depends polynomially on \ $s$.$\hfill \blacksquare$\\

  \begin{thm}\label{alpha inv.}
  Let \ $\lambda_1, \dots, \lambda_k$ \ be the fundamental invariants of a rank $k$ \ $[\lambda]-$primitive fresco such that \ $p_j \geq 1$ \ for each \ $j \in [1,k-1]$, and let \ $\mathcal{F}_0(\lambda_1, \dots, \lambda_k)$ \ the subset of \ $\mathcal{F}(\lambda_1, \dots, \lambda_k)$ \ of isomorphism classes of \ $E$ \  satisfying the condition that \ $F_{k-1}$ \ and \ $E\big/F_1$ \ are semi-simple, where \ $(F_j)_{j \in [1,k]}$ \ is the principal J-H. sequence of \ $E$. Then there exists a unique  polynomial function
  $$ \alpha : \mathcal{F}_0(\lambda_1, \dots, \lambda_k) \to \C  $$
  with the following properties :
  \begin{enumerate}[i)]
  \item \ $\alpha([E]) = 0 $ \ if and only if \ $E$ \ is semi-simple.
  \item When \ $E$ \ has a generator \ $e$ \ such that \ $(a - \lambda_{k-1}.b).(a - \lambda_k.b).e$ \ is in \ $F_{k-2}$, we have \ $\alpha(E) = \alpha(G_{k-1})$ \ where \ $G_{k-1} : = \A.(a - \lambda_{k-1}-1).b).e$.
  \item When \ $E \simeq \A\big/\A.(a - \lambda_1.b).(1 + x.b^{p_1})^{-1}.(a - \lambda_2.b) $ \ then \ $\alpha([E]) = x $.
  \end{enumerate}
  \end{thm}

\parag{proof} We define \ $\alpha$ \ for the rank \ $2$ \ case as the coefficient of \ $b^{p_1}$ \ in \ $S_1$ \ in any standard presentation \ $E \simeq \A\big/\A.(a - \lambda_1.b).S_1^{-1}.(a - \lambda_2.b)$. This is given by the easy computation which gives the rank \ $2$ \ classification for \ $[\lambda]-$primitive frescos,  and fulfills condition iii). Then thanks to the previous proposition we may define \ $\alpha(E)$ \ for any rank \ $k \geq 3$ \ $[\lambda]-$primitive fresco \ $E$ \ such that \ $F_{k-1}$ \ and \ $E\big/F_1$ \ are semi-simple via condition ii). The condition i) and the fact that \ $\alpha$ \  is a  polynomial function are proved in the induction step given by the proposition. $\hfill \blacksquare$\\

Let me give in rank \ $3$ \ a polynomial \ $F \in \C[s]$ \ such that we have \ $F(s) = \alpha([E(s)])$ \ for those \ $s$ \ for which \ $E(s)$ \ is in \ $\mathcal{F}_0(\lambda_1,\lambda_2,\lambda_3)$.\\
Assume \ $p_1\geq 1 $ and \ $p_2 \geq 1$. Then the necessary and sufficient condition to be in \ $\mathcal{F}_0(\lambda_1,\lambda_2,\lambda_3)$ \ is \ $s^1_{p_1} = s^2_{p_2} = 0 $.

\begin{lemma}
Let \ $E : = \A\big/\A.(a - \lambda_1.b).S_1^{-1}.(a - \lambda_2.b).S_2^{-1}.(a - \lambda_3.b)$ \ such that  \ $s^1_{p_1} = s^2_{p_2} = 0 $;  the complex number \ $ \alpha(E) $ \ is the coefficient of \ $b^{p_1+p_2}$ \ in \ $p_2.V.S_1$ \ for \ $j = 1,2$ \ where \ $V \in \C[[b]]$ \ is a solution of the differential equation 
 $$ b.V' = p_2.(V - S_2) .$$
\end{lemma}

The proof is left to the reader as an exercice.\\

This gives the following formula for \ $F$ :
$$ F(s) =  p_2.\sum_{j\not= p_2, j=0}^{p_1+p_2} \ s^1_{p_1+p_2-j}.\frac{s^2_{j}}{p_2-j}.  $$
Remark that for \ $S_2 = 1$ \ we find that \ $F(s)$ \ reduces to \ $s^1_{p_1+p_2}$. Using the commutation relation \ $(a - \lambda_2.b).(a - \lambda_3.b) = (a - (\lambda_3+1).b).(a - (\lambda_2-1).b)$ \ we find in this case that \ $E(s)$ \ has a normal rank \ $2$ \ sub-module with fundamental invariants \ $\lambda_1, \lambda_3+1$ \ and parameter \ $s^1_{p_1+p_2}$ \ which is precisely \ $\alpha(E(s))$.

\parag{Consequence} 
 If we have rank \ $4$ \ $[\lambda]-$primitive fresco \ $E$ \ with its principal J-H. sequence  such that \ $F_{j+1}\big/F_{j-1}$ \ is semi-simple for \ $j \in [1,3]$,  we may look inductively to \ $\alpha(F_3)$. If it is zero  then we may look at \ $\alpha(E/F_1)$, and if it is zero we may look at \ $\alpha(E)$. If this is zero then \ $E$ \ is semi-simple, and so,  the inductive  vanishing of all \ $\alpha-$invariants of sub-quotients of the principal J-H. sequence give a necessary and sufficient condition for semi-simplicity of \ $E$.\\
This  easily extends to any rank \ $k$ \ $[\lambda]-$primitive fresco. $\hfill \square$\\


Now we shall show that in the previous situation when the \ $\alpha-$invariant of \ $E$ is not zero it is the parameter of any normal rank \ $2$ \ sub-theme of \ $E$. Although such a rank \ $2$ \ normal sub-theme is not unique (in general), its isomorphism class is uniquely determined from the fundamental invariants of \ $E$ and from the \ $\alpha-$invariant of \ $E$.

\begin{prop}\label{sub-themes}
Let \ $E$ \ be a rank \ $k \geq 2$ \ $[\lambda]-$primitive fresco such that \ $F_{k-1}$ \ and \ $E\big/F_1$ \ are semi-simple and with \ $\alpha(E) \not= 0$. Note \ $\lambda_1, \dots, \lambda_k$ \ the fundamental invariants of \ $E$ \ and \ $p(E) : = \sum_{j=1}^{p-1} \ p_j$. Then there exists at least one rank \ $2$ \ normal  sub-theme in \ $E$ \ and  each rank \ $2$ \ normal sub-theme of \ $E$ \ is  isomorphic to
\begin{equation*}
 \A\big/\A.(a - \lambda_1.b).(1 + \alpha(E).b^{p(E)})^{-1}.(a - (\lambda_k+k-2).b). \tag{@}
 \end{equation*}
\end{prop}

\parag{proof} The case \ $k =2$ \ is clear, so we may assume that \ $k \geq 3$ \ and we shall prove the proposition by induction on \ $k$. So assume that the proposition is knwon for the rank \ $k-1 \geq 2$ \ and let \ $E$ \ be a rank \ $k$ \ $[\lambda]-$primitive fresco satisfying our assumptions. The fact that there exists a rank \ $2$ \ normal sub-theme is consequence of the induction hypothesis as  \ $G^{\tau}$ \ for any \ $\tau \in \C$ \ is normal rank \ $k-1$ \ in \ $E$ \ and satisfies again our assumptions.\\
Recall that the fundamental invariants of \ $G^{\tau}$ \ are \ $\mu_1, \dots , \mu_{k-1}$ \ with \ $\mu_j = \lambda_j$ \ for \ $j \in [1,k-2]$ \ and \ $\mu_{k-1} = \lambda_k + 1$. We have also \ $\alpha(G^{\tau}) = \alpha(E)$ \ for each \ $\tau$, thanks to the proof of the theorem \ref{alpha inv.}.  As we have  \ $p(G^{\tau}) = p(E)$ \ for each \ $\tau$, and \ $\mu_{k-1} = \mu_1 + p(E) - 1$ \ the inductive hypothesis  implies that any rank \ $2$ \ normal sub-theme of any \ $G^{\tau}$ \ is isomorphic to \ $(@)$. Then, to complete the proof, it is enough to show that any rank \ $2$ \ normal sub-theme of \ $E$ \ is contained in some \ $G^{\tau}$.
\smallskip

Let \ $T$ \ be a rank \ $2$ \ normal sub-theme of \ $E$. We shall first prove that its fundamental invariants are equal to \ $\lambda_1, \lambda_k+k-2$. As \ $E\big/F_1$ \ is semi-simple, the image of \ $T$ \ by the quotient map \ $E \to E\big/F_1$ \ has rank \ $1$ \ and this implies that \ $F_1\cap T$ \ is a rank \ $1$ \ normal sub-module. Then this implies that \ $F_1$ \ is the unique normal rank \ $1$ \ submodule of \ $T$, proving that \ $\lambda_1(T) = \lambda_1$.
\smallskip
We shall prove now that \ $\lambda_2(T) = \lambda_k+k-2$.\\
First remark that the uniqueness of the principal J-H. sequence of \ $E$ \ implies the uniqueness of the quotient map \ $E \to E\big/F_{k-1} \simeq E_{\lambda_k}$ \ because any surjective map \ $q : E \to E_{\lambda_k}$ \ will produce principal  a J-H. sequence for \ $E$ \ by adjoining a principal J-H. sequence for \ $Ker\, q$. So a quotient \ $E\big/H$ \ admits a surjective map on \ $E_{\lambda_k}$ \ if (and only if) \ $H \subset F_{k-1}$. So \ $H$ \ has to be semi-simple. Then the quotient \ $E\big/T$ \ has no surjective map on \ $E_{\lambda_k}$.\\
The exact sequence of frescos
$$ 0 \to T\big/F_1 \to E\big/F_1 \to E\big/T \to 0  $$
gives the equality  in \ $\A$ :
 $$P_{T\big/F_1}.P_{E\big/T} = P_{ E\big/F_1} = (a - \lambda_2.b) \dots (a - \lambda_k.b) .$$
 But as \ $E\big/T$ \ is semi-simple and does not has a quotient isomorphic to \ $E_{\lambda_k}$ \ this implies \ $T\big/F_1 \simeq E_{\lambda_k+ k-2}$, proving our claim.\\
 So \ $T$ \ is isomorphic to \ $\A\big/\A.(a - \lambda_1.b).(1 + \beta.b^{p(E)})^{-1}.(a - (\lambda_k+k-2).b)$ \ for some \ $\beta \in \C^*$. Let \ $x \in T$ \ be a generator of \ $T$ \ which is annihilated by 
$$(a - \lambda_1.b).(1 + \beta.b^{p(E)})^{-1}.(a - (\lambda_k+k-2).b).$$
 It satisfies 
\begin{equation*}
(a - (\lambda_k+k-2).b).x \in F_1 \subset F_{k-2} .\tag{*}
\end{equation*}
We shall determine all \ $x \in E$ \ which satisfies the condition \ $(*)$ \ modulo \ $F_{k-2}$. Fix a generator \ $e : = e_k$ \ of \ $E$ \ such that \ $e_{k-1} : = (a - \lambda_k.b).e_k$ \ is in \ $F_{k-1}$ \ and satisfies \ $(a - \lambda_{k-1}.b).e_{k-1} \in F_{k-2}$. Then we look for \ $U, V \in \C[[b]]$ \ such that
$$ x : = U.e_k + V.e_{k-1} \quad {\rm satisfies} \quad (a - (\lambda_k+k-2).b).x \in F_{k-2} .$$
This leads to the equations
\begin{align*}
& b^.U' - (k-2).b.U = 0 \\
& U + b^2.V' - (p_{k-1}+k - 3).V = 0 
\end{align*}
and so we get
\begin{equation*}
U = \rho.b^{k-2} \quad {\rm and} \quad V = \frac{\rho}{p_{k-1}}.b^{k-3} + \sigma.b^{p_{k-1}+k-3}
 \end{equation*}
 Note that for \ $\rho = 0$ \ we would have \ $x \in F_{k-1}$ \ and this is not possible for the generator of a rank \ $2$ \ theme as \ $F_{k-1}$ \ is semi-simple. Then assuming \ $\rho \not= 0$ \ we may write
 $$ x = \frac{\rho.b^{k-3}}{p_{k-1}}.\left[ p_{k-1}.b.e_k + e_{k-1} + \frac{p_{k-1}.\sigma}{\rho}.b^{p_{k-1}}.e_{k-1}\right] \quad {\rm modulo} \ \  F_{k-2} .$$
 Now recall that the generator of \ $G_{\tau}$ \ is given by
 $$ (a - (\lambda_{k-1}-1).b).\varepsilon(\tau) $$
 where \ $\varepsilon(\tau) = e_k + \tau.b^{p_{k-1}-1}.e_{k-1} $. A simple computation gives
 $$  (a - (\lambda_{k-1}-1).b).\varepsilon(\tau) = p_{k-1}.b.e_k + e_{k-1} + p_{k-1}.\tau.b^{p_{k-1}}.e_{k-1} \quad {\rm modulo} \ F_{k-2} .$$
 As we know that each \ $G_{\tau}$ \ contains \ $F_{k-2}$, we conclude that any such \ $x$ \ is in  \ $G_{\tau}$ \ for \ $\tau = \sigma\big/\rho $, concluding the proof. $\hfill \blacksquare$

Using duality we may deduce from this result the following corollary.

\begin{cor}\label{quot. rk 2}
In the situation of the previous proposition \ $E$ \ as a rank \ $2$ \ quotient theme and any such rank \ $2$ \  quotient theme is isomorphic to
$$ \A\big/\A.(a - (\lambda_1-k+2).b).(1 + \beta(E).b^{p(E)})^{-1}.(a - \lambda_k.b) $$
where
$$ \beta(E) : = (-1)^k\frac{p_1.(p_1+p_2)\dots (p_1+ \dots+p_{k-2})}{p_{k-1}.(p_{k-2}+p_{k-1})\dots (p_2+ \dots p_{k-1})} \alpha(E)$$
\end{cor}

The proof is left as an exercice for the reader who may use the following remark.

\parag{Remark} Let \ $E$ \ be a rank \ $2$ \ \ $[\lambda]-$primitive theme with fundamental invariants \ $\lambda_1, \lambda_2$ \ and parameter \ $\alpha(E)$. For \ $\delta \in \mathbb{N}, \delta \gg 1$ \ $E^*\otimes E_{\delta}$ \ is a \ $2$ \ $[1-\lambda]-$primitive theme with fundamental invariants\ $(\delta - \lambda_2, \delta - \lambda_1$ \ and parameter \ $-\alpha(E)$.\\
Then in the situation of the proposition \ref{sub-themes} \ $E^*\otimes E_{\delta}, \delta \gg 1$ \ satisfies again our assumptions and we have \ $\alpha(E^*\otimes E_{\delta}) = \beta(E)$.

  \section{The existence theorem for frescos.}
  The aim of this section  is to prove the the following existence theorem for  the fresco associated to a relative de Rham cohomology class :
   
   \begin{thm}\label{Existence}
 Let \ $X$ \ be a connected complex manifold of dimension \ $n + 1$ \ where \ $ n$ \ is a natural integer, and let \ $f : X \to D$ \ be an non constant proper  holomorphic function on an open  disc \ $D$ \ in \ $\mathbb{C}$ \ with center \ $0$. Let us assume that \ $df$ \ is nowhere vanishing outside of \ $X_0 : = f^{-1}(0)$.\\
 Let \ $\omega$ \ be a \ $\mathscr{C}^{\infty}-(p+1)-$differential form on \ $X$ \ such that \ $d\omega = 0 = df\wedge\omega $. Denote by \ $E$ \ the geometric \ $\A-$module \ $\mathbb{H}^{p+1}(X, (\hat{K}^{\bullet}, d^{\bullet}))$ \ and \ $[\omega]$ \ the image of \ $\omega$ \ in \ $E\big/B(E)$. 
 Then \ $\A.[\omega] \subset E\big/B(E)$ \ is a fresco.
 \end{thm}
 
 Note that this result is an obvious consequence of the finiteness theorem \ref{Finitude} that we shall prove below. It gives the fact that \ $E$ \ is naturally an \ $\A-$module which is of finite type over the subalgebra \ $\C[[b]]$ \ of \ $\A$, and so its \ $b-$torsion \ $B(E)$ \ is a finite dimensional \ $\C-$vector space. Moreover, the finiteness theorem asserts that \ $E \big/B(E)$ \ is a geometric (a,b)-module.

 \subsection{Preliminaries.}
 
 Here we shall complete and precise the results of the section 2 of [B.II]. The situation we shall consider is the following : let \ $X$ \ be a connected  complex manifold of dimension \ $n +1$ \ and \ $f : X \to \mathbb{C}$ \ a non constant holomorphic function such that \ $\{ x \in X / \ df = 0 \} \subset f^{-1}(0)$. We introduce the following complexes of sheaves supported by \ $X_0 : = f^{-1}(0)$
 \begin{enumerate}
 \item  The formal completion ''in \ $f$'' \ $(\hat{\Omega}^{\bullet}, d^{\bullet})$ \ of the usual holomorphic de Rham complex of \ $X$.
 \item The sub-complexes \ $(\hat{K}^{\bullet}, d^{\bullet})$ \ and \ $(\hat{I}^{\bullet}, d^{\bullet})$ \ of \  $(\hat{\Omega}^{\bullet}, d^{\bullet})$ \  where the subsheaves \ $\hat{K}^p$ \ and \ $\hat{I}^{p+1}$ \ are defined for each \ $p \in \mathbb{N}$ \  respectively as the kernel and the image of the map
 $$  \wedge df : \hat{\Omega}^p \to \hat{\Omega}^{p+1} $$
 given par exterior multiplication by \ $df$.  We have the exact sequence
 \begin{equation*} 0 \to (\hat{K}^{\bullet}, d^{\bullet}) \to (\hat{\Omega}^{\bullet},  d^{\bullet}) \to (\hat{I}^{\bullet}, d^{\bullet})[+1]  \to 0. \tag{1}
 \end{equation*}
 Note that \ $\hat{K}^0$ \ and \ $\hat{I}^0$ \ are zero by definition.
 \item The natural inclusions \ $\hat{I}^p \subset \hat{K}^p$ \ for all \ $p \geq 0$ \ are compatible with the diff\'erential \ $d$. This leads to an exact sequence of complexes
 \begin{equation*}
 0 \to (\hat{I}^{\bullet}, d^{\bullet}) \to (\hat{K}^{\bullet}, d^{\bullet}) \to ([\hat{K}/\hat{I}]^{\bullet}, d^{\bullet}) \to 0 .\tag{2}
 \end{equation*}
 \item We have a natural inclusion \ $f^*(\hat{\Omega}_{\mathbb{C}}^1) \subset \hat{K}^1\cap Ker\, d$, and this gives a sub-complex (with zero differential) of \ $(\hat{K}^{\bullet}, d^{\bullet})$. As in [B.07], we shall consider also the complex \ $(\tilde{K}^{\bullet}, d^{\bullet})$ \ quotient. So we have the exact sequence
 \begin{equation*}
  0 \to f^*(\hat{\Omega}_{\mathbb{C}}^1) \to (\hat{K}^{\bullet}, d^{\bullet}) \to (\tilde{K}^{\bullet}, d^{\bullet}) \to 0 . \tag{3}
  \end{equation*}
 We do not make the assumption here that \ $f = 0 $ \ is a reduced equation of \ $X_0$, and we do not assume that \ $n \geq 2$, so the cohomology sheaf in degree 1  of the complex \ $(\hat{K}^{\bullet}, d^{\bullet})$, which is equal to \ $\hat{K}^1 \cap Ker\, d$ \ does not coincide, in general with \ $f^*(\hat{\Omega}_{\mathbb{C}}^1)$. So the complex \ $ (\tilde{K}^{\bullet}, d^{\bullet})$ \ may have a non zero cohomology sheaf in degree 1.
  \end{enumerate} 
  Recall now that we have on the cohomology sheaves of the following complexes \\
   $(\hat{K}^{\bullet}, d^{\bullet}), (\hat{I}^{\bullet}, d^{\bullet}), ([\hat{K}/\hat{I}]^{\bullet}, d^{\bullet}) $ \ and \ $f^*(\hat{\Omega}_{\mathbb{C}}^1), (\tilde{K}^{\bullet}, d^{\bullet})$ \ natural operations \ $a$ \ and \ $b$ \ with the relation \ $a.b - b.a = b^2$. They are defined in a na{\"i}ve way by 
  $$  a : = \times f \quad {\rm and} \quad  b : = \wedge df \circ d^{-1} .$$
  The definition of \ $a$ \ makes sens obviously. Let me precise the definition of \ $b$ \ first in  the case of \ $\mathcal{H}^p(\hat{K}^{\bullet}, d^{\bullet})$ \ with \ $p \geq 2$  : if \ $x \in \hat{K}^p \cap Ker\, d$ \ write \ $x = d\xi$ \ with \ $\xi \in \hat{\Omega}^{p-1}$ \ and let \ $b[x] : = [df\wedge \xi]$. The reader will check easily that this makes sens.\\
  For \ $p = 1$ \ we shall choose \ $\xi \in \hat{\Omega}^0$ \ with the extra condition  that \ $\xi = 0$ \ on the smooth part of  \ $X_0$ \ (set theoretically). This is possible because the condition \\ 
  $df \wedge d\xi = 0 $ \ allows such a choice : near a smooth point of \ $X_0$ \ we can choose coordinnates such \ $ f = x_0^k$ \ and the condition on \ $\xi$ \ means independance of \ $x_1, \cdots, x_n$. Then \ $\xi$ \ has to be (set theoretically) locally constant on \ $X_0$ \ which is locally connected. So we may kill the value of such a \ $\xi$ \ along \ $X_0$.\\
  The case of the complex \ $(\hat{I}^{\bullet}, d^{\bullet})$ \ will be reduced to the previous one using the next  lemma.
  
  \begin{lemma}\label{tilde b}
  For each \ $p \geq 0$ \ there is a natural injective map
  $$\tilde{b} :  \mathcal{H}^p(\hat{K}^{\bullet}, d^{\bullet}) \to \mathcal{H}^p(\hat{I}^{\bullet}, d^{\bullet})$$
  which satisfies the relation \ $a.\tilde{b} = \tilde{b}.(b + a) $. For \ $p \not= 1$ \ this map is bijective.
  \end{lemma}
  
  \parag{Proof} Let \ $x \in \hat{K}^p \cap Ker\, d $ \ and write \ $x = d\xi $ \ where \ $x \in \hat{\Omega}^{p-1}$ \ (with \ $\xi = 0$ \ on \ $X_0$ \ if \ $p = 1$), and set \ $\tilde{b}([x]) : = [df\wedge \xi] \in \mathcal{H}^p(\hat{I}^{\bullet}, d^{\bullet})$. This is independant on the choice of \ $\xi$ \ because, for \ $p \geq 2$, adding \ $d\eta$ \ to \ $\xi$ \ does not modify the result as \ $[df\wedge d\eta] = 0 $. For \ $p =1$ \ remark that our choice of \ $\xi$ \ is unique.\\
   This is also independant of the the choice of \ $x $ \ in \ $ [x] \in \mathcal{H}^p(\hat{K}^{\bullet}, d^{\bullet})$ \  because adding \ $\theta \in \hat{K}^{p-1}$ \ to \ $\xi$ \ does not change \ $[df \wedge \xi]$.\\
   Assume \ $\tilde{b}([x]) = 0 $ \ in \ $ \mathcal{H}^p(\hat{I}^{\bullet}, d^{\bullet})$; this means that we may find \ $\alpha \in \hat{\Omega}^{p-2}$ \ such \ $df \wedge \xi = df \wedge d\alpha$. But then, \ $\xi - d\alpha $ \ lies in \ $\hat{K}^{p-1}$ \ and \ $x = d(\xi - d\alpha ) $ \ shows that \ $[x] = 0$. So \ $\tilde{b}$ \ is injective.\\
  Assume now \ $p \geq 2$.  If \ $df\wedge \eta $ \ is in \ $\hat{I}^p \cap Ker\, d$, then \ $df \wedge d\eta = 0 $ \ and \ $y : = d\eta $ \ lies in \ $\hat{K}^p \cap Ker\, d$ \ and defines a class \ $[y] \in  \mathcal{H}^p(\hat{K}^{\bullet}, d^{\bullet}) $ \ whose image by \ $\tilde{b}$ \ is \ $[df\wedge \eta] $. This shows the surjectivity of \ $\tilde{b}$ \ for \ $p \geq 2$.\\
   For \ $p=1$ \ the map \ $\tilde{b}$ \ is not surjective (see the remark below).\\
  To finish the proof let us  to compute \ $\tilde{b}(a[x] + b[x])$. Writing again \ $x = d\xi$, we get
   $$ a[x] + b[x] =[ f.d\xi + df \wedge \xi] = [d(f.\xi)] $$
   and so
   $$ \tilde{b}( a[x] + b[x] ) = [ df \wedge f.\xi ] = a.\tilde{b}([x]) $$
   which concludes the proof. $\hfill \blacksquare$
   
   \bigskip
   
   Denote by \ $i :  (\hat{I}^{\bullet}, d^{\bullet}) \to (\hat{K}^{\bullet}, d^{\bullet})$ \ the natural inclusion and define the action of \ $b$ \ on \ $\mathcal{H}^p(\hat{I}^{\bullet}, d^{\bullet})$ \ by \ $b : = \tilde{b}\circ \mathcal{H}^p(i) $. As \ $i$ \ is \ $a-$linear, we deduce the relation \ $a.b - b.a = b^2$ \ on \ $\mathcal{H}^p(\hat{I}^{\bullet}, d^{\bullet})$ \ from the relation of the previous lemma. \\
   
   The action of \ $a$ \ on the complex \ $ ([\hat{K}/\hat{I}]^{\bullet}, d^{\bullet}) $ \ is obvious and the action of \ $b$ \ is zero.\\
   
   The action of \ $a$ \ and \ $b$ \ on \ $f^*(\hat{\Omega}_{\mathbb{C}}^1) \simeq E_1\otimes \mathbb{C}_{X_0}$ \ are the obvious one, where \ $E_1$ \ is the rank 1 (a,b)-module with generator \ $e_1$ \ satisfying \ $a.e_1 = b.e_1$ \ (or, equivalentely, \ $E_1 : = \mathbb{C}[[z]]$ \ with \ $a : = \times z,\quad b : = \int_0^z $ \ and \ $e_1 : = 1$). \\
   Remark that the natural inclusion \ $f^*(\hat{\Omega}^1_{\mathbb{C}}) \hookrightarrow (\hat{K}^{\bullet}, d^{\bullet})$ \ is compatible with the actions of \ $a$ \ and \ $b$. The actions of \ $a$ \ and \ $b$ \ on \ $\mathcal{H}^1(\tilde{K}^{\bullet}, d^{\bullet}) $ \ are simply  induced by the corresponding actions on \ $\mathcal{H}^1(\hat{K}^{\bullet}, d^{\bullet})$.
   
   \parag{Remark} The exact sequence of complexes (1) induces  for any \ $p \geq 2$ \  a bijection
   $$ \partial^p : \mathcal{H}^p(\hat{I}^{\bullet}, d^{\bullet}) \to \mathcal{H}^p(\hat{K}^{\bullet}, d^{\bullet})$$
   and a short exact sequence 
   \begin{equation*}
    0 \to \mathbb{C}_{X_0} \to \mathcal{H}^1(\hat{I}^{\bullet}, d^{\bullet}) \overset{\partial^1}{\to} \mathcal{H}^1(\hat{K}^{\bullet}, d^{\bullet}) \to 0 \tag{@}
    \end{equation*}
   because of the de Rham lemma. Let us check  that for \ $p \geq 2$ \ we have \ $\partial^p = (\tilde{b})^{-1}$ \ and that for \ $p =1$ \ we have \ $\partial^1\circ \tilde{b} = Id$. If \ $x = d\xi \in \hat{K}^p \cap Ker\, d$ \ then \ $\tilde{b}([x]) = [df\wedge \xi]$ \ and \ $\partial^p[df\wedge\xi] = [d\xi]$. So \ $ \partial^p\circ\tilde{b} = Id \quad \forall p \geq 0$. For \ $p \geq 2$ \ and \ $df\wedge\alpha \in \hat{I}^p \cap Ker\, d$ \ we have \ $\partial^p[df\wedge\alpha] = [d\alpha]$ \ and \ $\tilde{b}[d\alpha] = [df\wedge\alpha]$, so \ $\tilde{b}\circ\partial^p = Id$. For \ $p = 1$ \ we have \ $\tilde{b}[d\alpha] = [df\wedge(\alpha - \alpha_0)]$ \ where \ $\alpha_0 \in \mathbb{C}$ \ is such that \ $\alpha_{\vert X_0} = \alpha_0$.  This shows that in degree 1 \ $\tilde{b}$ \ gives a canonical splitting of the exact sequence \ $(@)$.
   
  \subsection{$\A-$structures.}
   
   Let us consider now the \ $\mathbb{C}-$algebra
   $$ \A : = \{ \sum_{\nu \geq  0} \quad P_{\nu}(a).b^{\nu} \}$$
   where \ $P_{\nu} \in \mathbb{C}[z]$, and the commutation relation \ $a.b - b.a = b^2$, assuming that left and right multiplications by \ $a$ \ are continuous for the \ $b-$adic topology of \ $\A$.
   
   Define the following complexes of sheaves of left  \ $\A-$modules on \ $X$ :
   \begin{align*}
   & ({\Omega'}^{\bullet}[[b]], D^{\bullet})  \quad {\rm and} \quad   ({\Omega''}^{\bullet}[[b]], D^{\bullet}) \quad {\rm where} \tag{4}\\
    & {\Omega'}^{\bullet}[[b]] : = \sum_{j=0}^{+\infty} b^j.\omega_j \quad {\rm with} \quad \omega_0 \in \hat{K}^p  \\
     & {\Omega''}^{\bullet}[[b]] : = \sum_{j=0}^{+\infty} b^j.\omega_j \quad {\rm with} \quad \omega_0 \in \hat{I}^p \\
     & D(\sum_{j=0}^{+\infty} b^j.\omega_j) = \sum_{j=0}^{+\infty} b^j.(d\omega_j - df\wedge \omega_{j+1}) \\
     & a.\sum_{j=0}^{+\infty} b^j.\omega_j  = \sum_{j=0}^{+\infty} b^j.(f.\omega_j + (j-1).\omega_{j-1}) \quad {\rm with \ the \ convention} \quad \omega_{-1} = 0 \\
     & b.\sum_{j=0}^{+\infty} b^j.\omega_j  = \sum_{j=1}^{+\infty} b^j.\omega_{j -1}
     \end{align*}
     It is easy to check that \ $D$ \ is \ $\A-$linear and that \ $D^2 = 0 $. We have a natural inclusion of complexes of left \ $\A-$modules
     $$\tilde{i} :  ({\Omega''}^{\bullet}[[b]], D^{\bullet}) \to ({\Omega'}^{\bullet}[[b]], D^{\bullet}).$$
     
     Remark that we have natural morphisms of complexes
       \begin{align*}
     & u :   (\hat{I}^{\bullet}, d^{\bullet}) \to  ({\Omega''}^{\bullet}[[b]], D^{\bullet}) \\
     & v :  (\hat{K}^{\bullet}, d^{\bullet}) \to  ({\Omega'}^{\bullet}[[b]], D^{\bullet})
    \end{align*}
    and that these morphisms are compatible with \ $i$. More precisely, this means that we have the commutative diagram of complexes
 $$   \xymatrix{ (\hat{I}^{\bullet}, d^{\bullet}) \ar[d]^i \ar[r]^u & ({\Omega''}^{\bullet}[[b]], D^{\bullet}) \ar[d]^{\tilde{i}} \\
     (\hat{K}^{\bullet}, d^{\bullet})  \ar[r]^v &  ({\Omega'}^{\bullet}[[b]], D^{\bullet}) } $$
     
     The following theorem is a variant of theorem 2.2.1. of [B.II].

     \begin{thm}\label{(a,b)-structures}
     Let \ $X$ \ be a connected  complex manifold of dimension \ $n +1$ \ and \ $f : X \to \mathbb{C}$ \ a non constant holomorphic function such that 
      $$\{ x \in X / \ df = 0 \} \subset f^{-1}(0).$$
      Then the morphisms of complexes \ $u$ \ and \ $v$ \ introduced above are quasi-isomorphisms. Moreover, the isomorphims that they induce on the cohomology sheaves of these complexes are compatible with the actions of \ $a$ \ and \ $b$.
  \end{thm}
  
  \smallskip
  
  This theorem builds a natural structure of left \ $\A-$modules on each of the complex \\ 
  $(\hat{K}^{\bullet}, d^{\bullet}), (\hat{I}^{\bullet}, d^{\bullet}), ([\hat{K}/\hat{I}]^{\bullet}, d^{\bullet}) $ \ and \ $f^*(\hat{\Omega}_{\mathbb{C}}^1), (\tilde{K}^{\bullet}, d^{\bullet})$ \ in the derived category of bounded complexes of sheaves of \ $\mathbb{C}-$vector spaces on \ $X$.\\
   Moreover the short exact sequences
  \begin{align*}
&  0 \to (\hat{I}^{\bullet}, d^{\bullet}) \to (\hat{K}^{\bullet}, d^{\bullet}) \to ([\hat{K}/\hat{I}]^{\bullet}, d^{\bullet}) \to 0 \\
&  0 \to f^*(\hat{\Omega}_{\mathbb{C}}^1) \to (\hat{K}^{\bullet}, d^{\bullet}), (\hat{I}^{\bullet}, d^{\bullet}) \to (\tilde{K}^{\bullet}, d^{\bullet}) \to 0
\end{align*}
are equivalent to short exact sequences of complexes of left \ $\A-$modules in the derived category.

  \parag{Proof} We have to prove that for any \ $p \geq 0$ \ the maps \ $\mathcal{H}^p(u)$ \ and \ $\mathcal{H}^p(v)$ \ are bijective and compatible with the actions of \ $a$ \ and \ $b$. The case of \ $\mathcal{H}^p(v)$ \ is handled (at least for \ $n \geq 2$ \ and \ $f$ \ reduced) in prop. 2.3.1. of [B.II].  To seek for  completeness and for the convenience of the reader  we shall treat here the case of \ $\mathcal{H}^p(u)$.\\
  First we shall prove the injectivity of \ $\mathcal{H}^p(u)$. Let \ $\alpha = df\wedge \beta \in \hat{I}^p \cap Ker\, d$ \ and assume that we can find \ $U = \sum_{j=0}^{+\infty} b^j.u_j \in {\Omega''}^{p-1}[[b]]$ \ with \ $\alpha = DU$. Then we have the following relations
  $$ u_0 = df \wedge \zeta,  \quad \alpha = du_0 - df \wedge u_1 \quad {\rm and} \quad du_j = df\wedge u_{j+1} \ \forall j \geq 1. $$
  For \ $j \geq 1$ \ we have \ $[du_j] = b[du_{j+1}]$ \ in \ $ \mathcal{H}^p(\hat{K}^{\bullet}, d^{\bullet})$;  using corollary 2.2. of [B.II] which gives the \ $b-$separation of \ $\mathcal{H}^p(\hat{K}^{\bullet}, d^{\bullet})$, this implies \ $[du_j] = 0, \forall j \geq 1$ \ in \ $ \mathcal{H}^p(\hat{K}^{\bullet}, d^{\bullet})$. For instance we can find \ $\beta_1 \in \hat{K}^{p-1}$ \ such that \ $du_1 = d\beta_1$. Now, by de Rham, we can write \ $u_1 = \beta_1 + d\xi_1$ \ for \ $p \geq 2$, where \ $\xi_1 \in \hat{\Omega}^{p-2}$. Then we conclude that
   \ $ \alpha = -df\wedge d(\xi_1 + \zeta) $ \ and \ $[\alpha] = 0$ \ in \ $\mathcal{H}^p(\hat{I}^{\bullet}, d^{\bullet})$.\\
   For \ $p = 1$ \ we have \ $u_1 = 0 $ \ and \ $[\alpha] = [-df\wedge d\xi_1] = 0$ \ in \ $\mathcal{H}^1(\hat{I}^{\bullet}, d^{\bullet})$.\\
   We shall show now that the image of \ $\mathcal{H}^p(u)$ \ is dense in \ $\mathcal{H}^p({\Omega''}^{\bullet}[[b]], D^{\bullet})$ \  for its  \ $b-$adic topology. Let \ $\Omega : = \sum_{j=0}^{+\infty} \ b^j.\omega_j \in {\Omega''}^{p}[[b]]$ \ such that \ $D\Omega = 0$. The following relations holds \ $ d\omega_j = df\wedge \omega_{j+1} \quad \forall j \geq 0 $ \ and \ $\omega_0 \in \hat{I}^p$. The corollary 2.2. of [B.II] again allows to find \ $\beta_j \in \hat{K}^{p-1}$ \ for any \ $j \geq 0$ \ such that \ $d\omega_j = d\beta_j$. Fix \ $N \in \mathbb{N}^*$. We have
   $$ D(\sum_{j=0}^N b^j.\omega_j) = b^N.d\omega_N = D(b^N.\beta_N) $$
   and \ $\Omega_N : = \sum_{j=0}^N b^j.\omega_j  - b^N.\beta_N $ \ is \ $D-$closed and in \ ${\Omega''}^{p}[[b]]$. And we have \ $\Omega - \Omega_N \in b^N.\mathcal{H}^p({\Omega''}^{\bullet}[[b]], D^{\bullet})$, so the sequence \ $(\Omega_N)_{N \geq 1}$ \ converges to \ $\Omega$ \ in \ $\mathcal{H}^p({\Omega''}^{\bullet}[[b]], D^{\bullet})$ \ for its \ $b-$adic topology. Let us show that each \ $\Omega_N$ \ is in the image of \ $\mathcal{H}^p(u)$.\\
   Write \ $\Omega_N : = \sum_{j=0}^N b^j.w_j $. The condition \ $D\Omega_N = 0$ \ implies \ $dw_N = 0$ \ and \ $dw_{N-1} = df\wedge w_N = 0$. If we write \ $w_N = dv_N$ \ we obtain \ $d(w_{N-1} + df\wedge v_N) = 0$ \ and  \ $\Omega_N - D(b^N.v_N) $ \ is of degree \ $N-1$ \ in \ $b$. For \ $N =1$ \ we are left with \ $w_0 + b.w_1 - (-df\wedge v_1 + b.dv_1) = w_0 + df\wedge v_1$ \ which is in \ $\hat{I}^p \cap Ker\, d$ \ because \ $dw_0 = df\wedge dv_1$.\\
   To conclude it is enough to know the following two facts
   \begin{enumerate}[i)]
   \item The fact  that \ $\mathcal{H}^p(\hat{I}^{\bullet}, d^{\bullet})$ \ is complete for its \ $b-$adic topology.
   \item The fact that \ $Im(\mathcal{H}^p(u)) \cap b^N.\mathcal{H}^p({\Omega''}^{\bullet}[[b]], D^{\bullet}) \subset Im(\mathcal{H}^p(u)\circ b^N) \quad \forall N \geq 1 $.
   \end{enumerate}
   Let us first conclude the proof of the surjectivity of \ $\mathcal{H}^p(u)$ \ assuming i) and ii).\\
   For any \ $[\Omega] \in \mathcal{H}^p({\Omega''}^{\bullet}[[b]], D^{\bullet})$ \ we know that there exists a sequence \ $(\alpha_N)_{N \geq 1}$ \ in \ $ \mathcal{H}^p(\hat{I}^{\bullet}, d^{\bullet})$ \ with \ $\Omega - \mathcal{H}^p(u)(\alpha_N) \in b^N.\mathcal{H}^p({\Omega''}^{\bullet}[[b]], D^{\bullet})$. Now the property ii) implies that we may choose the sequence \ $(\alpha_N)_{N \geq 1} $ \ such that \ $[\alpha_{N+1}] - [\alpha_N]$ \ lies in \ $ b^N.\mathcal{H}^p(\hat{I}^{\bullet}, d^{\bullet})$. So the property i) implies that the Cauchy sequence \ $([\alpha_N])_{N \geq 1} $ \ converges to \ $[\alpha] \in \mathcal{H}^p(\hat{I}^{\bullet}, d^{\bullet})$. Then the continuity of \ $\mathcal{H}^p(u)$ \ for the \ $b-$adic topologies coming from its \ $b-$linearity, implies \ $\mathcal{H}^p(u)([\alpha]) = [\Omega]$.\\
   The compatibility with \ $a$ \ and \ $b$ \ of the maps \ $\mathcal{H}^p(u)$ \ and \ $\mathcal{H}^p(v)$ \ is an easy exercice. 
   
   \smallskip
   
   Let us now prove properties i) and ii).\\
   The property i) is a direct consequence of the completion of \  $\mathcal{H}^p(\hat{K}^{\bullet}, d^{\bullet})$ \ for its \ $b-$adic topology given by the corollary 2.2. of [B.II] \ and the \ $b-$linear isomorphism \ $\tilde{b} $ \ between \ $\mathcal{H}^p(\hat{K}^{\bullet}, d^{\bullet})$ \ and \ $\mathcal{H}^p(\hat{I}^{\bullet}, d^{\bullet})$ \ constructed in the lemma 2.1.1. above.\\
   To prove ii) let \ $\alpha \in \hat{I}^p \cap Ker\, d $ \ and \ $N \geq 1$ \ such that
   $$ \alpha = b^N.\Omega + DU $$
   where \ $\Omega \in {\Omega''}^{p}[[b]]$ \ satisfies \ $D\Omega = 0$ \ and where \ $U \in  {\Omega''}^{p-1}[[b]]$. With obvious notations we have
   \begin{align*}
   & \alpha = du_0 -df\wedge u_1\\
   & \cdots \\
   & 0  = du_j - df\wedge u_{j+1}  \quad \forall j \in [1, N-1] \\
   & \cdots \\
   & 0 = \omega_0 + du_N - df\wedge u_{N+1} 
   \end{align*}
   which implies \ $D(u_0+ b.u_1+ \cdots + b^N.u_N) = \alpha + b^N.du_N$ \ and the fact that \ $du_N$ \ lies in \ $\hat{I}^p \cap Ker \, d$. So we conclude that \ $[\alpha] + b^N.[du_N] $ \ is in the kernel of \ $\mathcal{H}^p(u)$ \ which is \ $0$. Then \ $[\alpha] \in b^N.\mathcal{H}^p(\hat{I}^{\bullet}, d^{\bullet})$.
   $\hfill \blacksquare$
   
     \parag{Remark} The map
  $$ \beta : ({\Omega'}[[b]]^{\bullet}, D^{\bullet}) \to ({\Omega''}[[b]]^{\bullet}, D^{\bullet})$$
  defined by \ $\beta(\Omega) = b.\Omega$ \ commutes to the differentials and with the action of \ $b$. It induces the isomorphism \ $\tilde{b}$ \ of the lemma \ref{tilde b} on the cohomology sheaves. So it is a quasi-isomorphism of complexes of \ $\mathbb{C}[[b]]-$modules.\\
 To prove this fact, it is enough to verify that the diagram\\
 $$ \xymatrix{&(\hat{K}^{\bullet}, d^{\bullet}) \ar[d]^{\tilde{b}} \ar[r]^v &  ({\Omega'}[[b]]^{\bullet}, D^{\bullet}) \ar[d]^{\beta} \\
 & (\hat{I}^{\bullet}, d^{\bullet}) \ar[r]^u &  ({\Omega''}[[b]]^{\bullet}, D^{\bullet})}$$
 
 \bigskip
 
 induces  commutative diagams on the cohomology sheaves. \\
  But this is clear because if \ $\alpha = d\xi$ \ lies in \ $ \hat{K}^p \cap Ker \, d$ \ we have \ $D(b.\xi) = b.d\xi - df\wedge \xi $ \ so \ $\mathcal{H}^p (\beta)\circ \mathcal{H}^p (v)([\alpha]) = \mathcal{H}^p (u)\circ \mathcal{H}^p (\tilde{b})([\alpha])$ \ in \ $\mathcal{H}^p ({\Omega''}[[b]]^{\bullet}, D^{\bullet}). \hfill \blacksquare$

    \subsection{The finiteness theorem.} 
  
  Let us recall some basic definitions on the left modules over the algebra \ $\A$.

  Now let \ $E$ \ be any left \ $\A-$module, and define \ $B(E)$ \ as the \ $b-$torsion of \ $E$. that is to say
  $$ B(E) : = \{ x \in E \ / \  \exists N \quad  b^N.x = 0 \}.$$
  Define \ $A(E)$ \ as the \ $a-$torsion of \ $E$ \ and 
   $$\hat{A}(E) : = \{x \in E \ / \ \mathbb{C}[[b]].x \subset A(E) \}.$$
   Remark that \ $B(E)$ \ and \ $\hat{A}(E)$ \ are  sub-$\A-$modules of \ $E$ \ but that \ $A(E)$ \ is not stable by \ $b$. 
  
  \begin{defn}\label{petit}
  A left \ $\A-$module \ $E$ \ is called {\bf small} when the following conditions hold
  \begin{enumerate}
  \item \ $E$ \ is a finite type \ $\mathbb{C}[[b]]-$module ;
  \item \ $B(E) \subset \hat{A}(E)$ ;
  \item \ $\exists N \ / \  a^N.\hat{A}(E) = 0 $ ;
  \end{enumerate}
  \end{defn}
  
  Recall that for \ $E$ \ small we have always the equality \ $ B(E) = \hat{A}(E)$ \ (see [B.I] lemme 2.1.2) and that this complex vector space is finite dimensional. The quotient \ $E/B(E)$ \ is an (a,b)-module called {\bf the associate (a,b)-module} to \ $E$.\\
  Conversely, any left \ $\A-$module \ $E$ \ such that \ $B(E)$ \ is a finite dimensional \ $\mathbb{C}-$vector space and such that \ $E/B(E)$ \ is an (a,b)-module is small.\\
  The following easy  criterium to be small will be used later :
  
  \begin{lemma}\label{crit. small}
  A left \ $\A-$module \ $E$ \ is small if and only if the following conditions hold :
  \begin{enumerate}
  \item \ $\exists N \ / \ a^N.\hat{A}(E) = 0 $ ;
  \item \ $B(E) \subset \hat{A}(E) $ ;
  \item  \ $\cap_{m\geq 0} b^m.E \subset \hat{A}(E) $ ;
  \item \ $Ker \, b$ \ and \ $Coker \, b$ \ are finite dimensional complex vector spaces.
  \end{enumerate}
  \end{lemma} 
  
  As the condition 3 in the previous lemma has been omitted in [B.II] (but this does not affect this article because this lemma was used only in a case were this condition 3 was satisfied, thanks to proposition 2.2.1. of {\it loc. cit.}), we shall give the (easy) proof.
  \parag{Proof} First the conditions 1 to 4  are obviously necessary. Conversely, assume that \ $E$ \ satisfies these four conditions. Then condition 2 implies that the action of \ $b$ \ on \ $\hat{A}(E)\big/B(E)$ \ is injective. But the condition 1 implies that \ $b^{2N} = 0$ \ on \ $\hat{A}(E) $ \ (see [B.I] ). So we conclude that \ $\hat{A}(E) = B(E) \subset Ker\, b^{2N}$ \ which is a finite dimensional complex vector space using condition 4 and an easy induction. Now \ $E/B(E)$ \ is a \ $\mathbb{C}[[b]]-$module which is separated for its \ $b-$adic topology. The finitness of \ $Coker \, b$ \ now shows that it is a free finite type \ $\mathbb{C}[[b]]-$module concluding the proof. $\hfill \blacksquare$
  
  \begin{defn}\label{geometric}
  We shall say that a left \ $\A-$module \ $E$ \ is {\bf geometric} when \ $E$ \ is small and when it associated (a,b)-module \ $E/B(E)$ \ is geometric.
  \end{defn}
  
  The main result of this section is the following theorem, which shows that the Gauss-Manin connection of a proper holomorphic function produces geometric \ $\A-$modules associated to vanishing cycles and nearby cycles.

 \begin{thm}\label{Finitude}
 Let \ $X$ \ be a connected complex manifold of dimension \ $n + 1$ \ where \ $ n$ \ is a natural integer, and let \ $f : X \to D$ \ be an non constant proper  holomorphic function on an open  disc \ $D$ \ in \ $\mathbb{C}$ \ with center \ $0$. Let us assume that \ $df$ \ is nowhere vanishing outside of \ $X_0 : = f^{-1}(0)$.\\
 Then the \ $\A-$modules 
 $$ \mathbb{H}^j(X, (\hat{K}^{\bullet}, d^{\bullet})) \quad {\rm and} \quad \mathbb{H}^j(X, (\hat{I}^{\bullet}, d^{\bullet})) $$
 are geometric for any \ $j \geq 0 $.
 \end{thm}
 
In the proof we shall use the \ $\mathscr{C}^{\infty}$ \ version of the complex \ $(\hat{K}^{\bullet}, d^{\bullet})$. We define \ $K_{\infty}^p$ \ as the kernel of \ $\wedge df : \mathscr{C}^{\infty,p} \to \mathscr{C}^{\infty,p+1}$ \ where \ $\mathscr{C}^{\infty,j}$ \ denote the sheaf of \ $\mathscr{C}^{\infty}-$ \ forms on \ $X$ \ of degree p, let \ $\hat{K}_{\infty}^p$ \ be the \ $f-$completion and \ $(\hat{K}_{\infty}^{\bullet}, d^{\bullet})$ \ the corresponding de Rham complex.

  The next lemma is proved in [B.II] (lemma  6.1.1.)
 
 \begin{lemma}\label{diff}
 The natural inclusion
 $$ (\hat{K}^{\bullet}, d^{\bullet}) \hookrightarrow (\hat{K}_{\infty}^{\bullet}, d^{\bullet}) $$
 induce a quasi-isomorphism.
 \end{lemma}
 
 \parag{Remark} As the sheaves \ $\hat{K}_{\infty}^{\bullet}$ \ are fine, we have a natural isomorphism
 $$ \mathbb{H}^p(X, (\hat{K}^{\bullet}, d^{\bullet})) \simeq H^p\big(\Gamma(X, \hat{K}_{\infty}^{\bullet}), d^{\bullet}\big).$$
 
 Let us denote by \ $X_1$ \ the generic fiber of \ $f$. Then \ $X_1$ \ is a smooth compact complex manifold of dimension \ $n$ \ and the restriction of \ $f$ \ to \ $f^{-1}(D^*)$ \ is a locally trivial \ $\mathscr{C}^{\infty}$ \ bundle with typical fiber \ $X_1$ \ on \ $D^* = D \setminus \{0\}$, if the disc \ $D$ \ is small enough around \ $0$. Fix now \ $\gamma \in H_p(X_1, \mathbb{C})$ \ and let \ $(\gamma_s)_{s \in D^*}$ \ the corresponding multivalued horizontal family of \ $p-$cycles \ $\gamma_s \in H_p(X_s, \mathbb{C})$. Then, for \ $\omega \in \Gamma(X, \hat{K}_{\infty}^p \cap Ker\, d)$, define the multivalued holomorphic function
 $$ F_{\omega}(s) : = \int_{\gamma_s} \frac{\omega}{df} .$$
 Let now 
  $$\Xi : = \oplus_{\alpha \in \mathbb{Q} \cap ]-1,0], j \in [0,n]} \quad  \mathbb{C}[[s]].s^{\alpha}.\frac{(Log s)^j}{j!} .$$
  This is an \ $\A-$modules with \ $a$ \ acting as multiplication by \ $s$ \ and \ $b$ \ as the primitive in \ $s$ \ without constant. Now if \ $\hat{F}_{\omega}$ \ is the asymptotic expansion at \ $0$ \ of \ $F_{\omega}$, it is an element  in \ $\Xi$, and we obtain in this way an \ $\A-$linear map
 $$ Int :   \mathbb{H}^p(X, (\hat{K}^{\bullet}, d^{\bullet})) \to H^p(X_1, \mathbb{C}) \otimes_{\mathbb{C}} \Xi .$$
 To simplify notations, let  \ $E : =  \mathbb{H}^p(X, (\hat{K}^{\bullet}, d^{\bullet}))$. Now using Grothendieck theorem [G.65], there exists \ $N \in \mathbb{N}$ \ such that \ $Int(\omega) \equiv 0 $, implies \ $a^N.[\omega] = 0$ \ in \ $E$.  As the converse is clear we conclude that \ $\hat{A}(E) =  Ker(Int)$. It is also clear that \ $B(E) \subset Ker(Int)$ \ because \ $\Xi$ \ has no \ $b-$torsion. So we conclude that \ $E$ \ satisfies properties 1 and 2 of the lemma \ref{crit. small}. The property 3 is also true because of the regularity of the Gauss-Manin connection of \ $f$.
 
 \parag{End of the proof of theorem \ref{Finitude}} To show that \ $E : =  \mathbb{H}^p(X, (\hat{K}^{\bullet}, d^{\bullet}))$ \ is small, it is enough to prove that \ $E$ \ satisfies the condition 4 of the lemma \ref{crit. small}. Consider now the long exact sequence of hypercohomology of the exact sequence of complexes
 $$   0 \to (\hat{I}^{\bullet}, d^{\bullet}) \to (\hat{K}^{\bullet}, d^{\bullet}) \to ([\hat{K}/\hat{I}]^{\bullet}, d^{\bullet}) \to 0 .$$
 It contains the exact sequence
 $$  \mathbb{H}^{p-1}(X, ([\hat{K}\big/\hat{I}]^{\bullet}, d^{\bullet})) \to \mathbb{H}^p(X, (\hat{I}^{\bullet}, d^{\bullet})) \overset{\mathbb{H}^p(i)}{\to} \mathbb{H}^p(X, (\hat{K}^{\bullet}, d^{\bullet})) \to \mathbb{H}^{p}(X, ([\hat{K}\big/\hat{I}]^{\bullet}, d^{\bullet})) $$
 and we know that \ $b$ \ is induced on the complex of \ $\A-$modules  quasi-isomorphic to \ $(\hat{K}^{\bullet}, d^{\bullet})$ \ by the composition \ $i\circ \tilde{b}$ \ where \ $\tilde{b}$ \ is a quasi-isomorphism of complexes of \ $\mathbb{C}[[b]]-$modules. This implies that the kernel  and the cokernel of \ $\mathbb{H}^p(i)$ \ are isomorphic (as \ $\mathbb{C}-$vector spaces) to \ $Ker\, b$ \ and \ $Coker \, b$ \ respectively. Now to prove that \ $E$ \ satisfies condition 4 of the lemma \ref{crit. small} it is enough to prove finite dimensionality for the vector spaces \ $  \mathbb{H}^{j}(X, ([\hat{K}\big/\hat{I}]^{\bullet}, d^{\bullet})) $ \ for all \ $j \geq 0 $.\\
 But the sheaves \ $[\hat{K}\big/\hat{I}]^j \simeq [Ker\,df\big/Im\, df]^j$ \ are coherent on \ $X$ \ and supported in \ $X_0$. The spectral sequence
 $$ E_2^{p,q} : = H^q\big( H^p(X, [\hat{K}\big/\hat{I}]^{\bullet}), d^{\bullet}\big) $$ 
 which converges to \ $ \mathbb{H}^{j}(X, ([\hat{K}\big/\hat{I}]^{\bullet}, d^{\bullet})) $, is a bounded complex of finite dimensional vector spaces by Cartan-Serre. This gives the desired finite dimensionality.\\
 To conclude the proof, we want to show that \ $E/B(E)$ \ is geometric. But this is an easy consequence of the regularity of the Gauss-Manin connexion of \ $f$ \ and of the Monodromy theorem, which are already incoded in the definition of \ $\Xi$ : the injectivity on \ $E/B(E)$ \ of the \ $\A-linear$ \ map \ $Int$ \ implies that \ $E/B(E)$ \ is geometric. \\
 Remark now that the piece of  exact sequence above gives also the fact that \ $\mathbb{H}^p(X, (\hat{I}^{\bullet}, d^{\bullet}))$ \ is geometric, because it is an exact sequence of \ $\A-$modules. $\hfill \blacksquare$
 
 \parag{Remark} It is easy to see that the properness assumption on \ $f$ \ is only used for two purposes : \\
 --To have a (global) \ $\mathscr{C}^{\infty}$ \ Milnor fibration on a small punctured disc around \ $0$, with a finite dimensional cohomology for the Milnor fiber.\\
 -- To have compactness of the singular set \ $\{df = 0 \}$, which contains the supports of the coherent sheaves \ $(Ker\, df\big/Im \,df)^i$.\\
 This allows to give with the same proof an analoguous finiteness result in many other situations.

  \section*{Bibliography}

\begin{itemize}

\item{[Br.70]} Brieskorn, E. {\it Die Monodromie der Isolierten Singularit{\"a}ten von Hyperfl{\"a}chen}, Manuscripta Math. 2 (1970), p. 103-161.

\item{[B.84]}  Barlet, D. \textit{Contribution du cup-produit de la fibre de Milnor aux p\^oles de \ $\vert f \vert^{2\lambda}$},  Ann. Inst. Fourier (Grenoble) t. 34, fasc. 4 (1984), p. 75-107.

\item{[B.93]} Barlet, D. {\it Th\'eorie des (a,b)-modules I}, in Complex Analysis and Geo-metry, Plenum Press, (1993), p. 1-43.

\item{[B.95]} Barlet, D. {\it Th\'eorie des (a,b)-modules II. Extensions}, in Complex Analysis and Geometry, Pitman Research Notes in Mathematics Series 366 Longman (1997), p. 19-59.

\item{[B.05]} Barlet, D. {\it Module de Brieskorn et forme hermitiennes pour une singularit\'e isol\'ee d'hypersuface}, revue de l'Inst. E. Cartan (Nancy) 18 (2005), p. 19-46. 

\item{[B.I]} Barlet, D. {\it Sur certaines singularit\'es non isol\'ees d'hypersurfaces I}, Bull. Soc. math. France 134 (2), ( 2006), p.173-200.

\item{[B.II]} Barlet, D. {\it Sur certaines singularit\'es d'hypersurfaces II}, J. Alg. Geom. 17 (2008), p. 199-254.

\item{[B.08]} Barlet, D. {\it Two finiteness theorem for regular (a,b)-modules}, preprint Institut E. Cartan (Nancy) (2008) $n^0 5$, p. 1-38, arXiv:0801.4320 (math. AG and math. CV)

\item{[B.09]} Barlet,D. {\it P\'eriodes \'evanescentes et (a,b)-modules monog\`enes}, Bollettino U.M.I. (9) II (2009), p. 651-697.

\item{[B.III]} Barlet, D. {\it Sur les fonctions \`a lieu singulier de dimension 1 }, Bull. Soc. math. France 137 (4), (2009),  p. 587-612.

\item{[B.10]} Barlet,D. {\it Le th\`eme d'une p\'eriode \'evanescente}, preprint Institut E. Cartan (Nancy) (2009) $n^0 33$, p. 1-57.

\item{[B.-S. 04]} Barlet, D. et Saito, M. {\it Brieskorn modules and Gauss-Manin systems for non isolated hypersurface singularities,} J. Lond. Math. Soc. (2) 76 (2007) $n^01$ \ p. 211-224.

\item{[Bj.93]} Bj{\"o}rk, J-E, {\it   Analytic D-modules and applications}, Kluwer Academic  publishers (1993).

\item{[G.65]} Grothendieck, A. {\it On the de Rham cohomologyof algebraic varieties} Publ. Math. IHES 29 (1966), p. 93-101.

 \item{[K.76]} Kashiwara, M. {\it b-function and holonomic systems}, Inv. Math. 38 (1976) p. 33-53.

\item{[M.74]} Malgrange, B. {\it Int\'egrale asymptotique et monodromie}, Ann. Sc. Ec. Norm. Sup. 7 (1974), p. 405-430.

\item{[M.75]} Malgrange, B. {\it Le polyn\^ome de Bernstein d'une singularit\'e isol\'ee}, in Lect. Notes in Math. 459, Springer (1975), p. 98-119.

\item{[S.89]} Saito, M. {\it On the structure of Brieskorn lattices}, Ann. Inst. Fourier 39 (1989), p. 27-72.

\end{itemize}

\end{document}